\documentclass[10pt]{article}

\usepackage{psfig}
\usepackage{epsfig}
\usepackage{color}

\usepackage{xypic,amsfonts,amsmath,amssymb}
\usepackage[english]{babel}

\author{Luc Pirio}
\title{Abelian functional equations, planar web geometry and
  polylogarithms\footnote{ 
{\small  This is a preliminary version. Any comments or remarks will
  be welcome.}         }}

\date{9 December 2002}
 
\newtheorem{theorem}{\textbf{Theorem}}

\newtheorem{definition}{\textbf{definition}}
\newtheorem{coro}{\textbf{corollary}}
\newtheorem{lemma}{\textbf{lemma}}
\newtheorem{conjecture}{\textbf {Conjecture}}
\newtheorem{prop}{\textbf{Proposition}}

\def \sm  {{ {\sf m} }}
\def \g {{\bf {\sf g} }}
\def \id { {\bf  I}_{\sf d}} 
\def \h {{\bf {\sf h} }}
\def \iv { {\bf  I}_{\sf v}} 
\def \l  #1 {{ {\bf L}{\mbox{i}}_#1}}

\def \eq #1 {( \cal {#1})} 
\def \rR {\mathbb R}

\def \rC {\mathbb C}
\def \rZ {\mathbb Z}
\def \rN {\mathbb N}
\def \tildhol #1 { \widetilde{{\cal O}} ( {#1} ) }
\def \hol #1 { {\cal O} ( #1 ) } 
\def \projc { \rC \mathbb P }
\def \llog { {\bf L}og }
\def \germhol #1  { \underline{ {\cal O}_ {#1} } }
\def \ll { {\bf L}}
\begin{document}

\maketitle



\begin{quote}
{\bf Abstract:} 
{In this paper we study abelian functional equations ({\sf Afe}),
  which are equations in the $F_i$'s of the type
$ F_1(U_1)+\dots+F_N(U_N)=0$. Here we restrict ourselves to the cases
when the $U_i$'s are rational functions in two variables. First we
prove that local measurable solutions actually are analytic and their
components are characterized as solutions to linear differential
equations constructed from the $U_i$'s. Then we propose two
``methods'' for solving ({\sf Afe}). Next we apply these methods to
the explicit resolution of generalized versions of classical
(inhomogeneous) {\sf Afe} satisfied by low order
polylogarithms. Interpreted in the framework of web geometry, these
results give us new non linearizable maximal rank planar webs
(confirming some results announced by G. Robert about one year ago).
Then we observe that there is a relation between these webs and
certain configurations of points in $ \projc^2$, which leads us to
define the notion of ``web associated to a configuration'': all these webs
seems to be of maximal rank. Finally, we apply the preceding results
to the problem of characterizing the dilogarithm and the trilogarithm
by the classical functional equation they respectively satisfy. In
particular, we show that, under weak regularity assumptions, the
trilogarithm is the only function which verifies the Spence-Kummer equation.}
\end{quote}

\section{Introduction and notations}

\subsection{Introduction}
In this paper, we undertake a general study of the general solutions 
$( F_1,..,F_N)$ of functional equations of the form
$$ \qquad   F_1(U_1(x,y))  +  F_2(U_2(x,y))  + \dots +
F_N(U_N(x,y))   =0     \qquad \qquad  \qquad \qquad  \qquad    ({ \cal E })       $$
where the $U_i$'s are real rational functions. We will call then
``abelian functional equations'' with real rational inner functions.
Such equations have appeared in mathematics a long time ago\,:
the equations 
\begin{align}
\qquad \qquad &  {\sf L}(\, x+y) ={  \sf L}(x)+  {  \sf L}  (y)
\quad  &&x,y \in \rR   \qquad  &&& ({\sf C}) \nonumber \\
\qquad \qquad   & \ll(\, xy)  =\ll(x)+\ll(y) \quad  &&  x,y>0 \qquad
&&& (C)  \nonumber 
\end{align} 
are respectively satisfied by any linear function and by the logarithm. From a historical point 
of view, equation $ ({C}) $ is closely related to the definition of
the logarithm itself and goes back to the 17th century. 

From the early 19th century onwards, many mathematicians  have gradually
discovered a particular class of special functions, the polylogarithms,
which verify some (inhomogeneous) functional equations of the 
type $
\eq{E} $ (see \cite{lewin}). Spence, Abel, Kummer (and others...) have
established numerous versions of the following functional equation
verified by the bilogarithm $ \l{2} $ for $ 0<x<y<1$ \,:
$$\ll ( x) -  \ll (y  ) - \ll ( \frac{x}{y}) -  \ll
 (\frac{1-y}{1-x})  +   \ll (
 \frac{x(1-y)}{y(1-x)} ) = -\frac{\pi ^2}{6} + \log(y)\:
 \log(\frac{1-y}{1-x}) \qquad 
    ({ L}_2)  $$
(this is Schaffer's form, see \cite{schaeff}). Spence and (mostly)
Kummer have discovered many functional equations satisfied by
polylogarithms of order less than 5, such as $ \l{3} $, which
verifies the following ``Spence-Kummer equation'', for $ 0<x<y<1 $, 
\begin{align}
  2 \, \ll (x) &  \, +  \, 2  \,\ll (y)  \, -   \,  \ll (\frac{x}{y})
  \,     
+ \, 2  \, \ll ( \frac{1-x}{1-y} ) \,   + \,  2 \ll (\frac{x(1-y)}{y(1-x)})    
- \, \ll  (xy) \nonumber \\
 &  \, +   \,  2 \ll (-  \frac{x(1-y)}{(1-x)}          ) 
+ \,  2 \,   \ll ( - \frac{(1-y)}{y(1-x)})   -  \, \ll (\frac{x(1-y)^2}{y(1-x)^2}) 
 \qquad   \quad   \qquad \qquad \qquad  ({ SK})  \nonumber \\
  &=2\l{3}(1)-\log(y)^2\log( \frac{1-y}{1-x})
+\frac{\pi^2}{3}\log(y)+\frac{1}{3}\log(y)^3   \nonumber
\end{align}
(we note $ {\sf R}_3(x,y) $ the right hand side of this equation).\\
The bilogarithm and some of its ``cousins'', such as the Rogers dilogarithm or the
single-valued Bloch-Wigner dilogarithm, are special functions which
have appeared in various branches of mathematics from the 1830's
onwards. Here are a few examples\,: in 1836, a result by Lobachevsky expressed the
volume of an ideal geodesic simplex in the 3-dimensional hyperbolic
space $ \mathbb H^3$ with vertices in $ \partial \mathbb H^3$ through
the dilogarithm. In
1935, G. Bol obtained the first example of a non-linearizable maximal rank planar 5-web by
considering the web associated to a (homogeneous) version of the
equation $ (L_2)$ of the bilogarithm (see part 4.1)\\
After a long period of neglect, for thirty years there has been an explosion of the occurrences of
polylogarithmic functions in many areas of mathematics (see, for
instance, \cite{gon0}, \cite{ost}, \cite{pol}, \cite{zagier1}, \cite{zagier}, 
 .. ).\\
Some mathematicians have generalized the construction of the
Bloch-Wigner dilogarithm to polylogarithms of any order: they have
constructed real univalued versions $ {\cal L}_n $ of $ \l{n} $,
defined and continuous on the whole $ \projc^1$. A theorem due to
Osterl{\'e}, Wojtkowiak and Zagier (Th{\'e}or{\`e}me 2 in  \cite{ost}) says that the $
 {\cal L}_{n} $'s verify
``clean'' versions of the functional equations satisfied by the
classical polylogarithms\,:
if for $ \l{n} $ we have an equation of the form
$$ \sum_{k=1}^N {\sf a}_k\, \l{n}(U_i)={\sf elem}_n $$
with ${\sf a}_k \in \mathbb C$, $U_k \in \rR(x,y)$, and $ 
{\sf elem}_n $ denoting a complex polynomial in some functions of the form $
{\bf L}_{i_{j_k}} \circ g_k $ with $1\leq j_k <n$ and $ g_k \in \rR(x,y)$
, then $ \sum_{k=1}^N {\sf a}_k\, ({\cal L }_{n}\circ U_i) $ is
constant.
This shows that the theory of the functional
equations of the polylogarithms can be considered a particular case
of the general study undertaken here.\\

Determining continuous
functions $ {\sf L}$ (resp. $ \ll $) satisfying $ ({\sf C}) $ (or
in an equivalent way $
({C}) $)  was important for early 19th century 
mathematicians\,: it allowed them to
justify in a rigorous manner the summation of ``Newtoon's binomial series'': $
{(1+x)}^{\alpha}= 1+\alpha x + \alpha( \alpha-1)x^2/2+...$ (with $ 
\alpha \in \rR \mbox{ and } |x|<1 $). Cauchy was the first to
 rigorously determine continuous solutions for these equations (now
known as ``Cauchy equations''). It was 
an application of the formalism that he had introduced into analysis
(see section 21.5 of \cite{aczdh}).
The problem of characterizing the  solutions of homogeneous versions
of the equations satisfied by polylogarithms is an interesting one.
Few results have been obtained in this direction (see part 4.2 and the conjecture
(1.6) in \cite{gangl})) although it could be an useful
way to retrieve certain results: see, for instance, the remark 4.1.2
in \cite{gelfandmac} or the proof of Lobachevsky's result stated above 
which is sketched in \cite{gonp} (page 7). \\

In this paper we study local solutions $ (F_1,..F_N) $ of the general equation $ \eq{E} $ at $
\omega \in \rR^2 $ and, in the spirit of the second part of Hilbert's
5th problem (see \cite{aczel}), we want
to make minimal assumptions of regularity on the $ F_i$'s so as to have
``nice properties'' for these solutions \,: measurability will appear 
natural (see 2.1). Under this assumption, we first prove (in
proposition 1 of part 2.2.1) that any
local solution of $ \eq{E} $ is in fact analytic (modulo a condition
of genericity on $\omega$): this allows us to complexify the problem
and to restrict ourselves to the study of local holomorphic solutions of the
complex version of $ \eq{E} $. As already noticed by Abel, one functional equation in several variables can
determine several unknown functions which must be very specific. In
our case, this ``philosophy'' works very well and gives us the \\

{\bf Theorem A} 
{\em Let be $ {\sf R}=(U_i) \in \rR(x,y)^N $ such that $ {\cal W}_{\sf
    R} $ is a web (i.e. the singular locus $\Sigma_{\sf R}  \subset
  \projc^2$ of $ {\sf R}$ is proper, see 1.2 for definitions). Let be
  $ \omega \in \rR^2 \setminus \Sigma_{\sf R}$ fixed. Then for each $
  i \in \{1,..,N\} $ there exists a linear differential equation 
$({\sf Lde}_i)$, the coefficients of which are algebraic
  functions such that if $F_1,..,F_N$ are measurable germs satisfying
  the equation $
F_1(U_1)+..+F_N(U_N)=0$ in a neighbourhood of $\omega $, then every
$F_i$ is analytic and generically satisfies  the equation $({\sf
  Lde}_i)$. The germ $F_i$ admits analytic continuation along any path
in the Zariski open set $X_i=U_i(\projc^2\setminus \Sigma) \subset
  \projc^1 $.}\\

Our result is explicit: for $ \Sigma_{\sf R} $, we have an explicit
formula in terms of the functions $U_i$'s. And, 
given a $N$-uplet $ {\sf R} $, we can explicitly construct the equation $ (
  {\sf Lde}_i)$ for every $i$ in terms of the $U_i$'s again.\\
 We prove this theorem by using
  mostly elementary methods of complex analysis. The proof can be
  divided into
 3 parts: from proposition 1, we know that the $F_i$'s are
  analytic germs. Then we complexify the setting. By successive
  differentiations along the level curves of the functions $U_i$'s, we
  construct for each i the linear differential equation $({\sf
    Lde}_i)$    from the equation $ \eq{E} $. This method is
  essentially an application to our case of Abel's method for solving
  functional equations in several variables, described in \cite{abel}. Finally we prove the analytic
  continuation along any path in $X_i$ by using a simple and general
  geometrical argument (see proposition~3).\\

 From the proof of this theorem, we deduce two methods to solve
equations of the form $ \eq{E} $. The first, called ``Abel's
method'', is explained in 2.3.1. It is effective and can be implemented
on a computer: it consists in solving the equation $ {\sf Lde _i}$
given by theorem A in order to reconstruct the solutions of $ \eq{E} $. The
second method, exposed in 2.3.2,  is not so general. It is based on
the fact that (modulo suitable condition on the $U_i$'s)
solutions of $ \eq{E} $ with logarithmic growth are characterized by
their monogromy, which can be determined~a~priori. 
In the third part, we first explicitly solve equations associated to
the classical equations of polylogarithms $(L_2)$ and $(SK)$ stated
above. Then in 3.5 we apply
Abel's method to an equation noted $ ( {\cal E}_{\sf c})$  associated
to a degenerate configuration $ {\sf c}$ of 5 points
in $\projc^2 $ (see figure 3). \\
In part 4.1, we interpret the preceding results in the
framework of planar web geometry: we obtain new ``exceptional webs''. 
In particular we prove the\\

{\bf Theorem B}
{\em The Spence-Kummer web $ {\cal W}_{\cal SK} $  associated
 to the equation $(SK)$ is an exceptional 9-web. }\\

The fact that we have found an explicit equivalent of the space of
abelian relations for this web in 3.4, allows us to study its sub-webs.
Thus we discover two non-equivalent exceptional 6-webs, and 
an exceptional 7-web. As in the case of Bol's web, numerous abelian relations for these
exceptional webs are constructed from polylogarithms. Then we observe
that, modulo a suitable change of coordinates, all these exceptional
webs are related to certain configurations of points in $\projc^2$ .

This remark leads us to define (see definition 4) the notion of ``web associated to a configuration of n points in the complex
projective plane''. 
Next we consider the web $ {\cal W}_{\sf c}$ associated to the
configuration  $  {\sf c}$. From the explicit basis of
solutions of $ ( {\cal E}_{\sf c})$ obtained in 3.5, we can now
construct a basis of the space of abelian relations of $ {\cal W}_{\sf
  c}$ showing that this web is exceptional. Then we state
some general results about webs associated to configurations of n points, for
$n=3,4,5$: \\

{\bf Theorem C}
{\em Let be $ n=3,4$ or $5$. The web associated to any (degenerate if $n=5$) configuration of $n$
points in $ \projc^2 $  is of maximal rank. Therefore it is exceptional if it
contains a sub-configuration of 4 points in general position. }\\

 This allows us to formulate a conjecture which could give numerous
 exceptional webs and therefore numerous equations of the form $ \eq{E} $.
Since the equations in part 3 (which are related to webs
associated to configurations) are mostly constructed by using
iterated integrals, this conjecture could give functional
equations for higher order polylogarithms.

 In part 4.2 we apply the preceding results to the problem of
 characterizing measurable functions $ \ll $ satisfying equations 
 $ (L_2)$ or $ (SK)$. 
We prove, with weak regularity assumptions, that $ \l{2} $ and 
$ \l{3} $  are characterized by these equations. In the case of the trilogarithm,
the result is new ( $ \l{3} $ is considered here as an analytic function on
 $ ]-\infty,1[\, $)\,:\\

{\bf Theorem D}
{\em Let $ F \! : \, ] -\infty,1 \: [ \: \rightarrow
\mathbb R $ be a measurable function such that for $ \, 
0<x<y<1 $, we have 
\begin{align}    
2  \, F(U_1(x,y)) & + \, 2 \, F(U_2(x,y))- \, F(U_3(x,y))      \nonumber \\
    & \qquad  +\,   2 \, F(U_4(x,y)) +\,   2 \, F(U_5(x,y)) -
  \, F(U_6(x,y)) \qquad    
\nonumber \\
 &   \qquad \qquad     + \, 2 \, F(U_7(x,y))+\,  2 \, F(U_8(x,y))-
   \,  F(U_9(x,y))= \, {\sf R}_3(x,y) \nonumber 
\end{align}
If $F$ is derivable at 0, then  $ F=\l{3} $ . }\\

This gives a proof of Goncharov's ``remark'' about the problem of
characterizing $ \l{3} $ by the Spence-Kummer equation, stated in
\cite{gon} (page 209). \\

{\bf remark :} {\bf 1.} This paper is an extended version of the
preprint \cite{pirio}. \\
{\bf 2.}While the author was working on the subject, he was told by G. Henkin 
that in a personal communication to him (nov. 2001), A. H{\'e}naut
announced that his colleague G. Robert had found that the Spence-Kummer's web
is of maximal rank by constructing an explicit basis of the space of abelian relations, which is equivalent
to part 3.4 of this paper. G. Robert had interpreted this in the framework of
web geometry and had obtained new exceptional d-webs for $d=6,7$
and $8$. But no additional information about this
has been given~until~now.\\
{\bf Acknowledgments\,:} The author would like to thank G. Henkin for
introducing him to this subject and discussing it with him. The
geometrical idea of the proof of proposition~3 comes from discussions with
J.M. Trepreau. Thank to C. Mourougane for his remarks and to A. Bruter
for her help to put this paper in form.

\subsection{Notations}
We introduce here some notations which we will use in the paper.\\

Throughout this paper, $N$ will be a fixed integer bigger than 3.\\

If no precision is given, for every $i=1,..,N$, $U_i$ will denote a
non-constant element of $ \rR(x,y)$ considered as a holomorphic map
 $ \projc^2\setminus S_i \rightarrow \projc^1$, where $S_i$ denotes the
 locus of indetermination of $U_i$\,: it is a finite set.\vspace{0.1cm}\\

A functional equation of the form $ F_1(U_1)+...+F_N(U_N)=0$ will be
called ``an abelian functional equation'' (ab. {\sf Afe}) with real
rational inner functions. The name comes from the notion of abelian
relation in web geometry, itself related to the notion of abelian
sum in algebraic geometry (see part 4.1.1 or part 2.2 in the expository paper \cite{hen1}).\\

In the whole text, $\eq{E} $ will denote a general ${\sf Afe}$ \ \ $ \sum_{i=1}^N
F_i(U_i)=0$.\vspace{0.1cm} \\

The foliation $ {\cal F}\{U_i\}$ (or more shortly $ {\cal F}_i$) will be
the global singular foliation of $\projc ^2$, the leaves of which are
the level curves of $U_i$. Let be ${\sf R}=(U_1,..,U_N) $ a $N$-uplet
of real rational functions. To the unordered set of foliations $ {\cal
  F}_{\sf R}=\{ {\cal F}_i \,
| \, i=1,..,N \} $, we associate the following algebraic subset of $ \projc^2$:
( ${\sf S}_i$ denotes the singular locus of the foliation $ {\cal F}_i$)
$$  \Sigma _{\sf R} := \Bigl( \cup_{i=1}^N  \, {\sf S}_{i}  \Bigr) \, \bigcup
\biggl( \cup_{ i \neq j } \left\{ \, \eta \in \projc ^2\setminus ( { S}_i\cup
{ S}_j)  \, | \,
( dU_i \wedge d   U_j)  \,      (\eta)=0 \, \right\} \biggr) $$
By definition, $ {\cal F}_{\sf R} $ is a web if $ \Sigma _{\sf R}$ is
proper in $ \projc^2$. In this case we note $ {\cal W}{ \{U_i \} }$
or $ {\cal W}_{{\sf R}}$ for  $ {\cal F}_{\sf R} $, and $ \Sigma
_{{\cal W}_{\sf R}} $ for $ \Sigma _{\sf R}$ and the latter will be
called the singular locus
of the web. Because $ \Sigma _{{\cal
  W}_{\sf R}}\,  $is the union of the singular locus of
the foliations ${\cal F}_i$  with the locus in
which the leaves of the foliations are not in general position, it
depends only on the web and not on the functions $U_i$.\\
The web $ {\cal W}\{ {\cal E} \}$ associated to $ \eq{E} $
will be the web  $ {\cal W}\{ U_i \}$. \\

If ${\cal F}$ is a sheaf of function germs on $\mathbb K \mathbb P
^d$ where $\mathbb K= \rR $ or $\rC$ and $ d=1,2$, $\underline{{\cal F}\scriptstyle{\omega}}$
will denote the function germs of this sheaf at $\omega \in
\mathbb K \mathbb P^d$ and we will note $ \underline{{{\cal
      F}}\scriptstyle{\omega}}( X)$ the space of determinations at
$\omega$ of the elements of  ${{{\cal F}}}(
X)$. In this paper, we will consider mostly the
sheaf ${\cal M}$ of measurable real valued function germs and the sheaf $ {\cal O}_X$ (ab. $ {\cal O}$) of holomorphic germs on a
complex manifold generally noted  $X$. Then $\widetilde{X}$ will be
the analytic universal covering of $X$, and $ \widetilde{ {\cal O}}_X $
(ab. $ \widetilde{ {\cal O}} $) will be the sheaf $ {\cal
  O}_{\widetilde{X}}$ of multivalued holomorphic functions on $X$. 
In the paper, $X$ will be a Zariski open set in $ \projc^k$ with $
k=1,2$. In part 2.2.2, we will use the sheaf of multivalued
holomorphic functions on $X$, with logarithmic growth at infinity,
noted $ \widetilde{{\cal O}}^{ \scriptscriptstyle{log}}_X  $ (ab.  $ \widetilde{{\cal O}}^{ \scriptscriptstyle{log}}  $).\\

If $ { \gamma}$ is a path linking $ \omega$ to $ \widetilde{\omega}$ in a
complex manifold $X$ and if $ {\sf K} \in \underline{ {\cal O}\scriptstyle{\omega}}$
admits an analytic continuation along $ {\bf \gamma}$, then we note $
{\sf K}^ { [{\bf \gamma}]}$ or $ {\cal M}_{ {\bf \gamma}} {\sf K}$ the
holomorphic germ at $ \widetilde{\omega}$ obtained by this analytic continuation.\\

If $ \omega \in \rR^2 $, then in the whole paper, we set $
\omega_i:=U_i(\omega) \in \rR \mathbb P^1 $ when it is well defined. Then a ``local solution of equation $
\eq{E} $ at $ \omega$ in the class ${\cal F}$'' will denote an element of the space 
$$ \underline{{\cal
      S} 
{ \stackrel{\cal F}{{}_\omega}}}   ( {\cal E})= \biggl\{
  (F_1,..F_N) \in \prod_{i=1}^N \underline{ {\cal F}\scriptstyle{\omega_i}} \;  |
  \;   \sum_{i=1}^N \,  F_i(U_i)=0 \; \mbox{ in } \, \underline{ {\cal
      F}\scriptstyle{\omega}} \; \biggr\} $$
We remark that if ${\cal F}={\cal O}$, then 
$ \underline{{\cal
      S} { \stackrel{\cal O}{{}_\omega}}}   ( {\cal E})$
is the space of the local
holomorphic solutions at $\omega$ of ``the complex version'' of $
\eq{E} $.\\
In the whole paper, to any  ${\bf H}=(H_1,..,H_N) \in \prod_i
\underline{ {\cal F}_{\omega_i}}  $ such that the sum $  \sum\,  F_i(U_i)$
is constant and equal to $c$, we associate the element $(H_1-c,..,H_N) $ of  $  
\underline{{\cal
      S} 
{ \stackrel{\cal F}{{}_\omega}}}   ( {\cal E})$ again noted~${\bf H} $.\\

If $J$ is a subset of $ \{1,..,N\}$, we note $ ( {\cal E}_J) $ the
equation $ \sum_{j \in J} F_j(U_j) $.\\
 For $ \omega \not\in \Sigma_{
  \cal E}$ we have $ \omega \not \in \Sigma_{ {\cal E}_J}$ and there
is a linear embedding 
$ \underline{{\cal
      S} \scriptstyle{
{ \stackrel{\cal F}{\omega}}}}   ( {\cal E}_J) \hookrightarrow 
\underline{{\cal
      S} \scriptstyle{
{ \stackrel{\cal F}{\omega}}}}   ( {\cal E})$. So we will consider
the local solutions of $ ( {\cal E}_J )$ as particular local solutions
of $ \eq{E} $. For $p \in \{3,..,N\} $ we note $ F^p \underline{{\cal
      S} \scriptstyle{{ \stackrel{\cal F}{\omega}}}}   ( {\cal E})$ the sum 
$ \sum  \underline{{\cal
      S} \scriptstyle{
{ \stackrel{\cal F}{\omega}}} }  ( {\cal E}_P) $ where $P$ runs 
 over all the subsets of $p$-elements in $ \{1,..,N\}$. An
element of $ F^p \underline{{\cal
      S} \scriptstyle{{ \stackrel{\cal F}{\omega}}} }   ( {\cal E}) \setminus 
F^{\scriptstyle{q}} \underline{{\cal
      S} \scriptstyle{{ \stackrel{\cal F}{\omega}}}}   ( {\cal E})$ (with
  $q=p-1$) will be
  called a solution of order $p$ of the equation $ \eq{E} $. A solution of
  order $ p<N$ will be called a ``sub-solution'', when a solution
  of order $N$ will be a ``genuine solution'' of $ \eq{E} $.\\
By definition ``a solution with logarithmic growth'' of $ \eq{E} $
will be an element of $\underline{{\cal
      S} \scriptstyle{
{ \stackrel{\widetilde{{\cal O}}}{\omega}}}}^{\scriptscriptstyle{log}} ( {\cal E})$.\\

The components of most known solutions of {\sf Afe} with
rational inner functions are constructed from iterated integrals. This
notion goes back to the work of K.T. Chen, in the 60's. We state here
the notations about iterated integrals used in the paper.\\

Let us note $X=\projc^2\setminus \Sigma_{ \cal W} $ and $
Z=\projc^1\setminus U_i(\Sigma_{\cal W})$ where $i$ is a fixed element
of $ \{1,..,N\}$.
There exists a finite number of distinct points $a_1,..,a_{M_i+1}$
in $\projc^1$ such that we have 
$Z=\projc^1\setminus\{ a_i \} $. We can always assume that
$a_{M_i+1}=\infty $ (we can substitute $ g \circ
U_i $ for $U_i$ with $g \in PGl_2(\rC)$ such that $ g(a_{M_i+1})=\infty $. This
doesn't change the nature of the problem.) We inductively define the
iterated integrals which are functions noted
$ {\bf L}_{x_{{i_1}}...x_{i_m}}$ with $i_k \in \{ 1,..,M_i\}$\,:
if $z\in Z$ and $ \gamma $ is a path in $Z$ from $ \omega_i$ to $z$
defining a point over $z$ in $\widetilde{Z}$ ,
then we set
$$  \qquad {\bf L}_{x_{i_0}x_{i_1}...x_{i_m}}(z, \gamma):= \int_{ \omega_i , \gamma}^z  \frac{ {\bf L}_{x_{i_1}x_{i_2}...x_{i_m}}(\xi)}{ a_{i_0}- \xi
} \; d\xi \quad , \quad  i_0,..,i_m \in \{1,..,M_i\}      $$
These functions are holomorphic functions on the analytic universal covering $
\widetilde{Z}$ of $Z$ .\\
We note $ {\cal I}_{ \{Z\} } $ (or $ {\cal I}_{ \{ a_i \} } $)  the
subspace of $ \widetilde{ {\cal O}}(Z)$ spanned by the constants and
the iterated integrals defined above. It is well defined: it doesn't
depend of the base point $ \omega_i$.\\

In part 3, we will use special notations for some elements of $ {\cal
  I}_{ \{ -1,0,1\} } $ that we describe now:\,
let be $\Omega := \rC \setminus ( \Delta_0 \cup \Delta_1 \cup
  \Delta_{-1}) $ where $\Delta _0$ , $\Delta _1$ and $ \Delta _{-1} $ 
are respectively the half-lines $ i \rR^-$,  $ 1+i \rR ^+$  and $ -1+i 
\rR ^-$  of $\rC$. Now $\Omega$ is simply connected and does not
  contain $0$ 
, $1$ and $-1$, so for any $z\in \Omega$ the value of any function
  defined by the expression below  is
  well defined if we integrate along
any path in $\Omega$:

\begin{align}
 {\bf L}_{x_0}(\bullet)& =\log(\bullet)  &{\bf L}_{x_1}(\bullet)
 &=-\log(1-\bullet) &
  {\bf L}_{x_{-1}}(\bullet)&=\log(1+\bullet)  \nonumber  \\
{\bf L}_{x_0x_1}(\bullet)& =\l{2} (\bullet)  &{\bf
 L}_{x_1x_0}(\bullet)&=\int_1^{\, \bullet} \frac{ 
{\bf L}_{x_0}(\zeta)}{1-\zeta} d\zeta &{\bf L}_{x_0x_{-1}}(\bullet)&=\int_0^{\, \bullet} \frac{ 
{\bf L}_{x_{-1}}(\zeta)}{\zeta} d\zeta  \nonumber \\
{\bf L}_{x_{-1}x_0}(\bullet)& =     \int_1^{\, \bullet} \frac{ 
{\bf L}_{x_0}(\zeta)}{1+\zeta} d\zeta            &{\bf
 L}_{x_1x_{-1}}(\bullet)&=\int_0^{\, \bullet} \frac{ 
{\bf L}_{x_{-1}}(\zeta)}{1-\zeta} d\zeta &{\bf L}_{x_{-1}x_{1}}(\bullet)&=\int_0^{\, \bullet} \frac{ 
{\bf L}_{x_{1}}(\zeta)}{1+\zeta} d\zeta  \nonumber \\
{\bf L}_{x_{\epsilon}x_{\epsilon}}(\bullet)&=\frac{1}{2}( {\bf
 L}_{x_{\epsilon}}(\bullet) )^2  \mbox{ for} \epsilon & = -1,0,1 &
\, &  \,  & \,  \nonumber
\end{align}
A polylogarithmic function will be a function constructed
from elements of~$ {\cal I}_{ \{0,1\} } $.

\section{General properties of the solutions of $ \eq{E} $ }
 \subsection{preliminary remarks} 

Our object is to study the solutions of an abelian functional equation
$\eq{E} $ with real rational inner functions. Using the notations introduced in the preceding part, we want to study
(and possibly determine) the space $   \underline{{\cal
      S}\scriptstyle{{ \stackrel{\cal F}{\omega}}}}   ( \cal E)   $ of local
  solutions of $ (\cal E)$ around $ \omega$ in the~class~$ \cal F $. What we want to prove is that, roughly speaking, the solutions of $ \eq{E} $
  are analytic, admit analytic continuation on a Zariski open set of $
\rC \mathbb P ^1$ and form a finite dimensional linear~space.\\

 But we have to make some  restrictions on $ \cal F $ and $
\omega $ to avoid pathological situations for the space $   \underline{{\cal
      S}{ \stackrel{\cal F}{{}_\omega}}}   ( \cal E)   $ : we have to
  deal with at least measurable functions and we have to take $\omega$
  outside of the singular locus $ \Sigma _{\cal W} $ of the web ${\cal W}\{ U_i
  \} $.\\

These two assumptions appear reasonable and quite natural if we
consider the following simple and classical examples: \\
First, let us consider the ``generalized Cauchy equation'' 
 $$ ( {\cal C})   \qquad 
  F_1(x)+F_2(y)+F_3(\frac{x}{y})=0$$ 
It is well known that the space of multiplicative functions $ F : \mathbb
R^{+ \star} \rightarrow \mathbb R $ is infinite dimensional. To
any such function corresponds a solution $(F,F,-F)$ of $  ( {\cal C} )
. $
 Such functions generally are not measurable: if the
function actually is, then it is constructed from the logarithm.
So, if no restriction on the regularity of the $F_i$'s is made, the space of
solutions can be infinite dimensional, which 
contrasts with the
measurable setting in which we have $ \mbox{dim}_{{}_{\scriptstyle{ \mathbb
    R}}}  \; \underline{{\cal
      S}\scriptstyle{{ \stackrel{\cal M}{\omega}}}}  ( {\cal C} )  = 3 \, $.
The assumption of measurability of the solutions of the general
equation $ \eq{E} $  appears natural.\\
According to this assumption  we can expect the solutions to have some good regularity properties
such as  analyticity : in our case, the ``only''
non constant measurable solution of $ \eq{C} $  is $ ( \log, \log,
-\log) $, which is  analytic indeed. But obtaining a precise local
version of this statement needs to make another assumption, about
$\omega$. If we take $\omega=0\in \rR ^2$ , then 
$ \mbox{dim}_{{}_{\scriptstyle{ \mathbb
    R}}}  \; \underline{{\cal
      S}\scriptstyle{{ \stackrel{ \mathbf C^0}{\omega}}}}  ( {\cal C}
  )  = 2 $ : actually $
  \eq{C} $  doesn't admit any non-constant analytic (and even
  continuous) solution at the origin. This comes from the fact that $ 0 $ belongs to the singular locus
of $ {\cal W}\{ x,y,xy\} $. Therefore the point $\omega$ must not belong to
this singular locus if we want it to have nice properties for the space $\underline{{\cal
      S}\scriptstyle{{ \stackrel{ \cal F}{\omega}}}}  ( {\cal C} ) $.\\

Another (more trivial) example of the pathologies which appear
if we don't make any assumption of genericity on $\omega$ is given by the
functional equation  $ G_1(x)+G_2(y)+G_3(x)=0 $ noted $( {\cal T})$. 
Here, the singular locus is the whole $ \rC \mathbb P ^2$ and the
local solutions of $( {\cal T})$ a priori don't admit any analytic continuation and form an infinite dimensional linear space. \\

These two elementary examples show that both the hypotheses of
measurability for the $F_i$'s and of genericity for $\omega$ are quite
natural and reasonable. 
From now on, we will always suppose that these hypotheses are satisfied.

\subsection{General properties of the measurable solutions of $ (\cal E) $ }
 
\subsubsection{Analyticity of the measurable solutions}
We prove now that any measurable local solution at a generic
point $\omega$ of $\eq{E} $ is in fact analytic.
\begin{prop}
Let be $\omega\in \mathbb R ^2 \setminus \Sigma $  and  $ {\bf
  F}=(F_1,..,F_N)\in  \underline{{\cal
      S}{ \stackrel{\cal M}{{}_\omega}}}   ( \cal E) $.
\\ Then each $F_i$ is in fact an analytic germ at $\omega_i$. Its 
complexification gives a germ $  F_i^{\scriptscriptstyle{\rC}}  \in \underline{
  {\cal O}_{\omega_i}}$ such that $ {\bf F}^{\scriptscriptstyle{\rC}}
:=(F_1^{\scriptscriptstyle{\rC}},..,F_N^{\scriptscriptstyle{\rC}})$ is a holomorphic
solution of $({\cal E})$~at~$\omega$~.
\end{prop}
{\bf proof:} 
By hypothesis we have $\omega \not \in \Sigma$, so it comes from
Theorem 3.3. of \cite{jarai} that the $F_i$'s are continuous germs at
 $\omega_i$. By elementary tools of
      integration it comes next that they are $ C^{\infty}$ smooth
      germs, 
so we have to prove that they are in fact analytic. \\
We obtain analyticity through the same method as J.L. Joly and J.Rauch 
in \cite{jolyrauch},
by formulating the equation $ \eq{E} $ in the form of a linear elliptic
$(N+2) \times N$ differential system.
      Then the analyticity of the $F_i$'s follows from classical
      results on the regularity of solutions of elliptic systems (see \cite{petro}). Finally the unicity principle
      implies that 
$ {\bf
        F}^{\scriptstyle{c}} \in    \underline{{\cal
      S}{ \stackrel{\cal O}{{}_\omega}}}   ( \cal E) $. $\blacksquare$ \\

{\bf remark:}  Under the assumption that the $F_i$'s are smooth enough, we
can show (see the next section) that each $F_i$ generically satisfies a linear
differential equation with analytic coefficients and so is analytic
by a classical result of ordinary differential equations. But this more
elementary way to prove analyticity is not useful because thus
it is not
easy to deal with the genericity condition.\\

So we have two $\mathbb R$-linear morphisms: \\
the first is just the restriction of the real
part of the holomorphic solutions of the complex version of $ \eq{E} $ to $ \mathbb R ^2 $ : 
$$\begin{array}{rrcl}
   {\bf \rho} : &   \underline{{\cal
      S}{ \stackrel{\cal O}{{}_\omega}}}   ( \cal E) &
{\longrightarrow} &  \underline{{\cal
      S}{ \stackrel{\cal M}{{}_\omega}}}   ( \cal E) 
  \nonumber \\ 
 \quad &  {\bf G} & \longmapsto &   \Re e( {\bf G}_{ | \scriptstyle{\mathbb
    R ^2} } )      
\end{array}$$
and the second is the complexification of the solutions given by proposition 1 \\
$$ \begin{array}{rrcl}
 {\bf \varrho} : &  \underline{{\cal
      S}{ \stackrel{\cal M}{{}_\omega}}}   ( \cal E) &
{\longrightarrow }& 
 \underline{{\cal
      S}\scriptstyle{ \stackrel{\cal O}{\omega}}}   ( \cal E) \nonumber \\ 
 \quad &  {\bf F} & \longmapsto &   {\bf F}
  ^{\scriptscriptstyle{\rC}} 
\end{array} $$

It is clear that $   {\bf \varrho } \circ  {\bf \rho} = {\bf I}{\sf d}_{
  \scriptstyle{\underline{ {\cal S}{}_{\omega}^{\cal M}}( {\cal E})
  }} \, $, so the study of the measurable solutions of $({\cal E})$ at $\omega$ amounts to the study of the holomorphic local solutions of $ \eq{E} $ .

\subsubsection{characterization of the components of the holomorphic
  solutions of $ \eq{E} $ }

It is well known that, in the generic case, there is no non-constant
holomorphic solution of a general abelian functional equation.\\
Let us consider now the very specific case when $ \eq{E} $ has a non
trivial local holomorphic solution ${\bf F}=(F_1,..,F_N)$. Any
non-constant component germ $F_i$ of ${\bf F}$ must be a 
function of a very specific kind. The point is that any germ $F_i$ is just a local determination of a globally
defined but ramified function which satisfies a linear differential
equation with algebraic~coefficients. We formulate this in the following
\begin{theorem} 
 Let $N \geq3$ be an integer and $ {\bf \sf R}=(U_1,..,U_N)\in
\rR(x,y)^N $ be such that $ \Sigma_{\bf \sf R}$ is proper. 
Let $\omega \in \rR^2\setminus \Sigma_{\bf \sf R}$ be fixed. 
Then for every $ i \in \{1,..,N\} $ there exists a linear differential
equation $ ({\sf Lde}_i) $, the coefficients of which are algebraic
functions (meromorphic in a neighbourhood of $ \omega_i$), such that for
all $ (F_1,..,F_N) \in 
\underline{{\cal
      S}\scriptstyle{ \stackrel{\cal O}{\omega}}}   ( \cal E) $, the
  germ $F_i$ satisfies  $ ({\sf Lde}_i) $ in a neighbourhood
  of $ \omega_i$. The germ $F_i$ is a local
  determination at $\omega_i$ of a globally defined multivalued
  function on $ \projc^1$, the ramification points of which belong to
  the finite set $U_i(\Sigma_{\bf \sf R})\subset \projc ^1$. 
\end{theorem}

{\bf proof:} Without any loss of generality, we can assume that
$\omega=(0,0) \not \in \Sigma_{\bf \sf R}$  and $U_i(\omega)=0$ for $
i=1,..,N$. Let be $N$ germs $F_i\in \underline{ {\cal O}}({\rC},0) $ such
that $ {\sum_1^N F_i(U_i)=0} $ in a neighbourhood of~$\omega$ . 
For $\rho >0$ let's note $ {\sf D}_{\rho}=\{ z \in \rC \, | \,
|z|<\rho \, \} $ . \\
If $i \neq j $, since $\omega \not \in \Sigma_{\sf {\bf R}}$, $
(U_i,U_j) $ defines a system of holomorphic coordinates on a
neighbourhood $ \Omega_{ij} $ of $ \omega $. It is clear that we can find $ \epsilon>0 $ such that each $F_i$ is
holomorphic on the whole ${\sf D}_{\epsilon  }$, and such that 
$ \Omega:= \bigcap_k U_k ^{-1}( {\sf D}_{\epsilon  }) \subset
\Omega_{ij}$ for all $i \neq j $ .\vspace{0.2cm} \\ 
We now want to deduce from the functional equation $ \eq{E} $ a
linear differential equation ({\sf Lde}) satisfied by $F_N$ ( or by any other $F_i$, the
process remaining the same). To do this, we will find it useful to
introduce a more general class of equations than {\sf Afe} :
\begin{definition}
{\em Let be $ (N,M_1,..,M_N) \in \rN^{\star}\! \times \rN^N$, and let 
$ V_i\, , {\cal A}_{ij}  \, ( 1\leq i \leq N ,0\leq j \leq M_i) $ be 
holomorphic functions on an open set $  \Theta \subset \rC^2 $ .
An ``Abelian Differential Functional Equation''(ab. {\sf Adfe}) is an
equation of 
the type
$$ 
 \qquad \qquad     \sum_{i=1}^N \sum_{j=0}^{M_i} {\cal A}_{ij} \, G_i ^{(j)}(V_i)=0 
\qquad \qquad ({\cal A}_{\sf dfe}) 
$$
where the unknowns are the function germs $G_1,..,G_N \, $  which are
supposed smooth enough ($ G_i^{(k)}$
denoting the $k$th-derivative of $G_i$ for $ k \in \rN $).\\
If ${\cal A}_{iM_i}\not \equiv 0 $ for all $i$'s, then the $N$-uplet
$(M_1,..,M_N)$ is called ``the true type'' of the equation $ ({\cal A}_{\sf
  dfe}) $ , and its ``type'' if not.}
\end{definition}    
The notion of  {\sf Adfe} generalizes {\sf Afe} and {\sf Lde}\,: {\sf
  Afe} are {\sf Adfe} of the type $(0,...,0)$ and {\sf
  Lde} are {\sf Adfe} of the type $(M_1)$ with $M_1>0$.\\

 Let us assume that for all  for $i\neq j$, the couple $(V_i\, ,\! V_j)$
 (noted after definition~1) defines holomorphic coordinates on
$\Theta$.
We now describe a process to obtain an {\sf Adfe} of the type
$(M_1-1,M_2+1,...,M_N+1)$, or $(M_2+1,...,M_N+1)$ if $ \, M_1=0$, 
from an {\sf Adfe} of the true type $(M_1,..,M_N)$.\\

To begin with, let's study the case when $M_1>0$.\\
By definition we have ${\cal A}_{1M_1} \not \equiv 0$ on $\Theta$, 
therefore the equation $ ({\cal A}_{\sf
  dfe}) $ implies that on $ \Theta'=\Theta \setminus \! \{ {\cal
  A}_{1M_1}=0 \} $, so we have 
$$  G_1^{(M_1)}(V_1)+\sum_{j=0}^{M_1-1} 
\frac{ {\cal A}_{1j} }{ {\cal A}_{1M_1}} \, G_1^{(j)}(V_1)  
+\sum_{i=2}^N \sum_{j=0}^{M_i} 
\frac{ {\cal A}_{ij} }{ {\cal A}_{1M_1}} \, G_i^{(j)}(V_i)  
$$
Let ${\partial}$ be the vector field on $\Theta$ which corresponds to
the differentiation with respect to $V_2$ in the coordinate system $(V_1,V_2)$.

By application of this derivation to this last form of $ ({\cal A}_{\sf
  dfe}) $ we get a new {\sf Adfe} on $ \Theta'$:
$$ \sum_{j=0}^{M_1-1} 
{\partial}  (\frac{ {\cal A}_{1j} }{ {\cal A}_{1M_1}})G_1^{(j)}(V_1)  
+\sum_{i=2}^N \sum_{j=0}^{M_i} \left( 
{\partial} (   \frac{ {\cal A}_{ij} }{ {\cal A}_{1M_1}})
\, G_i^{(j)}(V_i)  
+ \frac{ {\cal A}_{ij} }{ {\cal A}_{1M_1}} {\partial}(V_i)  G_i^{(j+1)}(V_i)
\right)=0 
$$
which can be written
$$ 
  \qquad \qquad    \sum_{i=1}^N \sum_{j=0}^{\widetilde{M_i}} \widetilde{{\cal
         A}_{ij}} \, G_i ^{(j)}(V_i)=0  \qquad \qquad ({\cal A}_{\sf
  dfe}^2) 
$$
where $\widetilde{M_i}=M_1-1 \mbox{ (resp. } M_i+1) 
\mbox{ if } i=1  \mbox{ (resp. } i>1) $ and 
$$ (\star) \qquad \widetilde{{\cal A}_{ij}} =
\begin{cases}   \; {\partial}  (\frac{ {\cal A}_{1j} }{ {\cal A}_{1M_1}}) 
\qquad \qquad \qquad \qquad \mbox{ if } i=1 \\ 
\; {\partial}  (\frac{ {\cal A}_{i0} }{ {\cal A}_{1M_1}}) \qquad \qquad
\qquad  \qquad \mbox{ if } i>1
\mbox{ and } j=0 \\ 
  \; {\partial}  (\frac{ {\cal A}_{ij} }{ {\cal A}_{1M_1}})+ \frac{ {\cal
    A}_{ij-1} }{ {\cal A}_{1M_1}}{\partial}(V_i) \qquad \; \mbox{ if } 1<i 
\mbox{ and } 0<j\leq M_i \\ 
\;  \frac{ {\cal
    A}_{iM_i} }{ {\cal A}_{1M_1}}{\partial}(V_i) \qquad \qquad \qquad
    \quad  \mbox{ if } 1<i 
\mbox{ and } j=M_i+1 
\end{cases}$$
We remark that, because $ {\partial}(V_i)\neq0$ for $i>1$, no $ 
\widetilde{{\cal A}}_{iM_i+1}$ is a null function, so the
equation that we obtain is of the true type $ (K,M_2+1,....,M_N+1) $,
$K$ being an integer smaller than $ M_1-1 $ .

If $M_1=0$, then we similarly get an equation of the form
$$
\qquad \qquad    \sum_{i=2}^N \sum_{j=0}^{\widehat{M_i}} \widehat{{\cal
         A}_{ij}} \, G_i ^{(j)}(V_i)=0  \qquad \qquad ({\cal A}_{\sf
  dfe}^2) 
$$

where $\widehat{M_i}=M_i+1 $ for $ 2\leq j \leq N$ and 
$$ (\star \star) \quad \widehat{{\cal A}_{ij}} =
\begin{cases} 
\; {\partial}  (\frac{ {\cal A}_{i0} }{ {\cal A}_{1M_1}}) \qquad \qquad
\qquad  \qquad \mbox{ if } i\geq2
\mbox{ and } j=0 \\ 
  \; {\partial}  (\frac{ {\cal A}_{ij} }{ {\cal A}_{1M_1}})+ \frac{ {\cal
    A}_{ij-1} }{ {\cal A}_{1M_1}}{\partial}(V_i) \qquad \; \mbox{ if }
2\leq i 
\mbox{ and } 0<j\leq M_i \\ 
\;  \frac{ {\cal
    A}_{iM_i} }{ {\cal A}_{1M_1}}{\partial}(V_i) \qquad \qquad \qquad
    \quad  \mbox{ if } 2\leq i 
\mbox{ and } j=M_i+1 
\end{cases}$$
As in the preceding case, it's quite obvious  that we obtain an {\sf
  Adfe} of the true type $ (M_2+1,....,M_N+1)$.\\

In both cases ($M_1=0)$ or $ ( M_1>0$), we can apply these
operations again to $ ({\cal A}_{\sf
  dfe}^2) $.\\
After several applications  of this process  on $\Omega$ to the {\sf
  Afe}  $ \eq{E} $  
 we obtain an {\sf Adfe} of the type $(K)$ (with $ K\in \rN^{\star} $) on $\Omega'=\Omega  \setminus
 \Lambda$, where $ \Lambda $ is an analytic subset of $\Omega$ . This
equation can be written in the coordinate system $(U,V)=(U_{N-1},U_N)$
in the following form: 
$$ { A}_1(U,V)F_N^{(1)}(V)+{ A}_2(U,V)F_N^{(2)}(V)+            ....+
{ A}_K(U,V)F_N^{(K)}(V)=0 \qquad ( {\cal A}_{\sf dfe}^N)$$
Let us take now $(U_0,V_0) \in \Omega'$. By fixing $U=U_0$ in the preceding
equation, we get, in a neighbourhood of $V_0$, a linear differential
equation of order $K$ in the variable $V$, the solutions of which contain
$F_N$:
$$ \qquad \quad {\sf A}_1(V)F_N^{(1)}(V)+{\sf A}_2(V)F_N^{(2)}(V)+            ....+
{\sf A}_K(V)F_N^{(K)}(V)=0 \qquad \qquad ({\sf L}{\sf de}_N) $$
It is clear that this equation  $ ({\sf L}{\sf de}_N) $ doesn't depend
on the solution $(F_1,..,F_N)$ but only on the $U_i$'s. Then the $N$th
component of every solution $ {\bf F} \in \underline{{\cal
      S}\scriptstyle{ \stackrel{\cal O}{\omega}}}   ( \cal E) $ will
  verify this equation, at least generically in a neighbourhood of
  $\omega_i $.\\

From now on, we assume that we can take $(U_0,V_0)=(0,0)$. 
We now prove by ``induction on the type'' that the coefficients of the
preceding equation are algebraic functions of $V$.
Let be $ \rC\{U,V\}^{\scriptstyle{alg}}={\{ h \in \rC\{U,V\}        | \, \exists
\, Q\in \rC[U,V,W] \, Q(U,V,h (U,V))=0 \} } $.\\
It is well known that this space has some strong properties of closure:
\begin{prop}Let be $F,G \in \rC\{U,V\}^{\scriptstyle{alg}}$.
Then $F+G$, $F\times G$, $
\partial_U F ,\partial_V F $ and $ 1/F$ (if $F(0,0)\neq 0$) are still
elements of $ \rC\{U,V\}^{\scriptstyle{alg}}$.
If $ \Phi=(F,G)$ defines a germ of diffeomorphism of $ \rC^2 $ at the origin,
then the components of the local inverse $ \Phi^{-1}$ are algebraic
functions too.
\end{prop}
Let us note ${\partial}^{kl}_k$ the derivation on $ \Omega $ with respect
to $U_k$ in the coordinate system $(U_k,U_l)$. One can easily prove
that we have $ {\partial}^{kl}_k={\sf U}_k^{kl} \, \partial_U+ {\sf
  V}_k^{kl} \, \partial_V $ with  $  {\sf U}_k^{kl}, {\sf V}_k^{kl} \in
\rC\{U,V\}^{\scriptstyle{alg}}$. Then by proposition 5, the latter is
closed under the action of the ${\partial}^{kl}_k$'s. By proposition 5
again and from the above relations $(\star)$ and $(\star \star)$, if the ${\cal A}_{ij} $'s    of $ {\cal A}_{\sf dfe}$ are
algebraic functions, then the $\widetilde{\cal A}_{ij} $'s ( or the
$\widehat{\cal A}_{ij} $'s) of $ ({\cal A}_{\sf dfe}^2) $
 are still algebraic. Because all the
coefficients of ${\sf Afe} $ $ \eq{E} $ are equal to $1$, we get, by
induction, that the $ { A}_i$'s of $({\cal A}_{\sf dfe}^N)$ are elements of 
$ \rC\{ U, V \}^{alg}$. Therefore the ${\sf A}_i$'s of $({\sf L}{\sf
  de}_N)$   are algebraic functions of $V$.\\

Because the ${\sf A}_i$'s are algebraic, they are globally
defined but ramified. A classical result of the theory of linear
differential equations of a complex variable implies that the germ $F_N$ can be analytically
extended along any curve in $ \projc ^1 \setminus  { R}$, where $R$
is the union of the poles with the ramification points of the~${\sf
  A}_i$'s.\\

But this argument didn't allow us to prove that $ F_N $ admits
analytic continuation along any path in the whole $U_N(\projc^2
\setminus \Sigma_{ {\sf R}} )$ because, if 
it's not hard to see that the ramification points of the $ {\sf
  A}_i$'s are in $ U_N(\Sigma_{ {\sf R}})$ , it is not the same for
their possible poles, which can generate some ramification
for any solution of~$({\sf L}{\sf de}_N)$ .\\
The last part of the theorem comes from the following proposition
3. $\blacksquare$ 

\begin{prop} Let $X$ be a connected paracompact complex manifold of
  dimension $2$ and let $U_i : \, X \rightarrow \rC , \,  (i=1,..,N) $ be holomorphic
  functions such that, if $  i\neq j  $, we have $ dU_i \wedge dU_j
  \neq 0 $ on $ X$. If for $\omega \in X $ we have $ N$ holomorphic
  germs $F_i$ such that $ \sum_1^N F_i\circ U_i $ is a holomorphic germ at
  $\omega$ which can be analytically continued along any path in the whole
  $X$, then every $ F_i$ can be analytically continued along any path in $
  U_i(X) $ .
\end{prop}
{\bf Proof :} we will prove this proposition under the assumption that 
$ \sum_1^N F_i\circ U_i =0$. The proof in the general case is
similar.\\
For $ i=1,..,N$, let's note $ \Psi_i:= F_i \circ U_i \in \underline{
  {\cal O}}{}{\scriptstyle{\omega_i}}$.
Because $X$ is supposed paracompact, it is metrisable as a topological
space. We fix a metric on X, compatible with its topology. Then there
  exists $ \epsilon >0 $ such that each $\Psi_i$ is defined on $ B(\omega,
\epsilon) \subset \subset  X$. First we prove the following
\begin{lemma}
Let us assume that $X$ is an open ball in $ \rC ^2 $ centered in $\omega$, of radius
$ \rho \geq  \epsilon $ .
Then each $\Psi_i$ can be analytically extended to  $ X=B_{\rho}:=B(\omega,\rho) $.
\end{lemma}
{\bf proof of the lemma: } Let be $ \tau:=\mbox{sup}\{ \, \delta \in [
\epsilon , \rho] \; | \;  \mbox{ each } \Psi_i \mbox{ extends to }
B_{\delta}  \; \} $. \\
We want to prove that $\tau = \rho$ .
Let us suppose that $ \tau < \rho $ : by definition each $ \Psi_i $
extends analytically to $ B_{\tau }$. We note again $ \Psi_i $ this
extension.

Let us choose arbitrarily $ \eta \in \partial B_{\tau}$. 
We are going to prove that all the $\Psi_i$'s have a holomorphic
extension in a neighbourhood of $ \eta$ . By compacity, it will imply
that each $\Psi_i$ extends to a neighbourhood of the closure $
\overline{B_{\tau}} $ , which will contradict the definition of $ \tau
$.\\

Let $ (x,y) $ denote the standard complex coordinates on $ \rC ^2 $
.\\
We introduce the holomorphic vector fields of differentiation along the
level curves of the $U_i$'s : $ {\cal X}_i:= ( \frac{{\partial}
U_i}{ \partial y} ) \, \partial _x  - ( \frac{{\partial}
U_i}{ \partial x} ) \, \partial _y $. \\

According to the definition of $ \Psi_i$, we have $ {\cal X}_i \Psi_i=0 $ on $
B_{\epsilon }$ and therefore on $ B_{\tau  }$  by unicity theorem : $\Psi_i $
is constant along the level curves of $ U_i$ in $ B_{\tau  }$ . But
these level curves are globally defined on $X$ and in particular in a
neighbourhood of $\eta$~.\\
 This fact combined with the general position
assumption on these level curves at $\eta$ (formulated by $ dU_i \wedge
dU_j (\eta) \neq 0 $  according to the hypothesis of the theorem) will allow us
to extend each $ \Psi_i$ near $\eta$.\\

But we have to make it more precise:\\
We note $ \mbox{T}_{\! \eta} \, \partial B_{\tau} $ the real tangent space
of $ \partial B_{\tau} $ in $\eta$. It is a real subspace of real dimension $3$
of  the complex tangent space to $ \rC ^2$ at $\eta \, $, noted $
\mbox{T}_{\eta} \rC ^2$. It contains an unique complex line noted  $ \mbox{T}_{\! \eta}^{\rC}
\, \partial B_{\tau} $.\\
Let us extend $\Psi_1$ in a neighbourhood of
$\eta$.

Let $  \mbox{C}^j_{\eta} $ be the level curve of $U_j$ through
$\eta$. Since $ \mbox{d}U_1(\eta) \neq 0$, we know that there exists a
neighbourhood $ {\cal V} $ of  $ \eta$ such that $  \mbox{C}^1_{\eta}  \cap
{\cal V} $ is a complex 1-dimensional manifold. Let $ \mbox{T}_{\eta} ^\rC \mbox{C}^1 $ be
its holomorphic tangent space in $ \eta$.\\
Let us assume that $ \mbox{T}_{\! \eta}^{\rC}
\, \partial B_{\tau} $ and $ \mbox{T}_{\eta} ^\rC \mbox{C}^1 $ are
transverse ( i.e. their intersection  in $ \mbox{T}_{\eta} \rC
^2$ is null). Because all geometrical objects considered here are 
analytic, therefore smooth,
this condition of transversality, called `` condition $ ({\cal T}) $ '', 
 is open: there exists an open connected neighbourhood $ \mbox{V}_{\eta} \subset X$ of
 $\eta$ such that for all $\zeta \in \mbox{V}_{\eta} \cap \, \partial B_{\tau}$
 , the transversality condition between $  \mbox{C}^1_{\zeta} $ and
 $\partial B_{\tau} $ remains  satisfied.\\

Let us note $ \mbox{W}_{\eta}:= \mbox{V}_{\eta}  \cap \, U_1^{-1}\left(
U_1(\mbox{V}_{\eta} \cap \, \partial B_{\tau}) \right) $. It is an open neighbourhood of $\eta $.\\
For $ \zeta \in  \mbox{W}_{\eta} $ ,  then $  \mbox{C}^1_{\zeta} \cap \,
\partial B_{\tau}\neq \emptyset $. On the  other hand, $
\mbox{C}^1_{\zeta} $ verifies the transversality
condition $ ({\cal T})$. The fact that $ \mbox{T}_{\! \eta} \, \partial
B_{\tau} $ contains an unique complex line implies, for dimensional
reasons, that $ \mbox{T}_{\! \eta} \, \partial B_{\tau} \cap \,
\mbox{T}_{\eta} ^\rC \mbox{C}^1 \neq (0) $ .\\
We deduce that  $  \mbox{C}^1_{\zeta} \cap \,  B_{\tau}\neq \emptyset $ .
Then let us consider $\zeta' \in \mbox{C}^1_{\zeta} \cap \,  B_{\tau} $ : we
define the value of $\Psi_1$ in $\zeta$ by setting $
\Psi_1(\zeta):=\Psi(\zeta') $. Because $\Psi_1$ is constant along the
level curves of $U_1$ in $ \mbox{W}_{\eta} \cap \, B_{\tau} $ , it
comes that $ \Psi_1(\zeta)$ is well defined.\\
We remark that we have $ {\cal X}_1 \Psi_1=0 $ near $\eta$ again 
for this extension, so we have holomorphically extended  $ \Psi_1$ to $
B_{\tau} \cup \mbox{W}_{\eta} $ .\\

Let us suppose now that the condition $ ({\cal T}) $ is not satisfied by
$\mbox{C}^1_{\eta} $ . \\
It means that $  \mbox{T}_{\! \eta} ^\rC \partial B_{\tau} = 
\mbox{T}_{\eta} ^\rC \mbox{C}^1 $ . But the hypothesis $ dU_1 \wedge
dU_j (\eta) \neq 0 $ for $ j\geq 2 $ has the geometrical
interpretation that the curves $\mbox{C}^1_{\eta} $ and
$\mbox{C}^j_{\eta} $ are transverse in $\eta$ ( for $j\geq 2$ ). 
Therefore all the level curves $\mbox{C}^j_{\eta} $ (for $j \geq 2$) satisfy the
transversality condition $ ({\cal T}) $ at $\eta$ . By the same
argument than above, we can extend analytically each $\Psi_j$
($j\geq2)$ to a neighbourhood $\mbox{W}$ of $\eta$.
To extend $\Psi_1$ close to $\eta$ we will set $ \Psi_1:= -\sum_{j=2}^N
\Psi_j $ on W, which will~do.$\blacksquare $   \hspace*{0.1cm}\hfill{\bf end of
  the lemma's proof} \\
Now let's prove proposition 3.\\
Let $ \gamma : \, [0,1] \rightarrow X_1:=U_1(X)$ be a path with
$\omega_1$ as its origin. We want to extend $\Psi_1$ along~$\gamma$.
Because $ dU_1 \, \wedge \, dU_2\neq 0$  on  $X$ ,  $\Gamma=(U_1,U_2)$ defines
holomorphic coordinates in a neighbourhood of any point of $X$ .
This implies first that $ dU_1 \neq 0 $  on  $X$, so we can find a lift
$ \tilde{\gamma}$ of $\gamma$ to~$X$ through $U_1$ with $\omega $ as
its origin.\\

Since the 
support $| \tilde{\gamma} | $ of $ \tilde{\gamma}$ is compact, it
comes too that
we can find a subdivision $ \alpha_{-1}< 0=\alpha_0<\alpha_1<...< \alpha_M=1 $ of
$[0,1]$ such that, for every $j \in \{ 0,..,M-1\}$, there exists a
holomorphic chart $ ( \Theta _j , \Gamma | _{\Theta_j} )$ centered at 
$ \tilde{\gamma}(\alpha_j)$, such that $ \tilde{\gamma}
([ \alpha_{j-1},  \alpha_{j+1}] ) \subset  \Theta_j$, and  
$ \Gamma(\Theta_j) $ is a ball in
$\rC ^2$ with $\Gamma( \alpha_j)$  as its center.
We can apply  the precedent lemma when taking
$X= \phi_0(\Theta_0) $ and considering  the functions $ U_i \circ
{\phi_0}^{-1} $ instead of the functions~$U_i$.
So it comes that the $\Psi_i$ can be extended along $
\tilde{\gamma}|_{[\alpha_0,\alpha_1]} $ .\\
By iterating this process
$M-1$ times, we  finally get an extension of each $\Psi_i$ along $
\tilde{\gamma} $, again noted $\Psi_i$. \\

This gives us the analytic extension of $ \Psi_1 $ along  $ \gamma $ :
for each chart $ ( \Theta _j , \Gamma | _{\Theta_j}   )$ we have on $
{\Theta_j}$ some holomorphic vector
fields $ {\cal X}_i ^j $ of differentiation along the
level curves of $U_i$ . There
is $ g_i ^j \in {\cal O} ^{\star}( \Theta_j \cap \, \Theta_{j+1} ) $
such that $  {\cal X}_i ^j = g_i ^j {\cal X}_i ^{j+1} $ on 
$ \Theta_j \cap \, \Theta_{j+1} $ for all $i$ and $j<M$ . \\
By the construction of the lemma , we have $
{\cal X}_1 \Psi_1 =0 $ on each $ \Theta_j$. The holomorphic
inverse function theorem implies that we can write 
$\Psi_1=F_1^j(U_1)$ on each $ \Theta_j$ , where $F_1 ^j $ is holomorphic in a
neighbourhood of $ ( U_1 \circ \tilde{\gamma}) ([ \alpha_{j-1},
\alpha_{j+1}])={\gamma}([ \alpha_{j-1} ,\alpha_{j+1}] ) $ .
It is not difficult now to see that on $  \gamma([\alpha_{j} , \alpha_{j+1} ]) $
we have $F_1 ^j = F_1 ^{j+1} $ ( for $0<j<M-2 $ ) . \\
Now $F_1^0$ is the extension to $ U_1(\Theta_0) $ of the
original germ $F_1$. By setting $F_1:= F_1^j $ on $ U_1( \Theta_j)$
for $j=1,..,M-1$ , we
get an analytic extension of $F_1$ along $\gamma$ .\\
By the same way we can construct such an analytic continuation for every
$F_j$. $\blacksquare$ \\ \hspace*{0.1cm}\hfill{\bf end of
  the proof of proposition 3} \\
{\bf remarks:}\\
{\bf 1.} From the preceding proof, we get, using the same notations:
\begin{coro}
In the generic case, there are no non-constant local holomorphic
solutions of $\eq{E} $ at any $ \omega \not \in \Sigma_{\bf \sf R}$.
\end{coro}{\bf proof:} Because equation $ ({\cal A}_{\sf dfe}^N)$ is
of the true type $(K)$ with $K>1$, we can assume that $ A_K \equiv
1$. Let us assume that one of the $A_i$'s ($i=1,..,K-1$) depends on the
variable $U$. Then by differentiating with respect to $U$, we reduce 
$ ({\cal A}_{\sf dfe}^N)$ to an equation $ ({{\cal A}_{\sf dfe}^N} ')$
of the true type $ (K')$ with $ 0\leq K'<K$. The obstacle to this process
of reduction is when $K'=0$: the differentiation with respect to $U$ gives a trivial equation. It corresponds
to the cancellation of all the $ \partial_U
A_i$. It corresponds to $K-1$ (non-linear) differential conditions on the
$U_i$'s. Then the possibility to reduce $ ({\cal A}_{\sf dfe}^N)$ to
an $ \sf Adfe $ of the form $ \widehat{A}_N(U,V) F_N'(V)=0$, with $
\widehat{A}_N \not \equiv 0$, corresponds to the non-vanishing of a finite
number of differential expressions in the $U_i$'s. An element $(U_i)
\in \rR(x,y)^N$ satisfying these conditions will be said N-generic, 
and generic if $(U_{ \sigma(i)})$ is N-generic for every $ \sigma \in
\mathfrak S _N$. It is clear that the genericity condition described
here implies that any holomorphic solution of $ \eq{E} $ is
constant. $\blacksquare$\\

{\bf 2.} Theorem 1 implies that $  \underline{{\cal
      S}{ \stackrel{\cal O}{{}_\omega}}}   ( \cal E) $  is finite
  dimensional. \\
The following proposition (due to G. Bol, see \cite{bol}) 
gives an effective bound to its dimension and will be important in part 4.
\begin{prop}
If $ \omega \not \in \Sigma_{\sf R}$ , then $ \dim_{\scriptstyle{\rC}} 
  \, \underline{{\cal
      S}{ \stackrel{\cal O}{{}_\omega}}}   ( {\cal E}) \leq
  {N(N-1)}/{2} $ and this majoration is optimal.
\end{prop}
The majoration is just a particular case (for rational inner functions) of 
one of the first basic results of web geometry (see
\cite{blabol}). If we consider the case when $U_i(x,y)=x-a_i \,y$ 
for $ i=1,..,N$ , with $ 0=a_1<a_2<..<a_N=1$, we easily see 
that the bound $N(N-1)/2$ is reached in this case.\\
{\bf 3.} Some of the preceding results remain valid in a more general
situation. For instance, if instead of taking the $U_i$'s rational, we consider
some analytic germs, then the $H_i$'s solutions of $ \sum
H_i(U_i)=0$ generically verify a linear differential equation which
can be constructed from the $U_i$'s, and we find again
that, in the generic case, there is no non-constant solution.\\
{\bf 4.} The method used here to obtain a linear differential equation from a
functional equation is the one described by Abel in his first publication
\cite{abel}. We will call it `` Abel's  Method''.\\
{\bf 5.} The point is that this method is effective: for a given
$N$-uplet of rational functions, we can explicitly find a linear
differential equation satisfied by any component of any local
solution. We can even do this in an algorithmic way: see the next section.\\
Similarly, for a fixed $N$, we can explicitly find sufficient
conditions on $(U_i)_{i\leq N}$ so that there is no non-constant solution
of $ \eq{E} $.\\
{\bf 6.} Through the process used to extend the $F_i$'s in
proposition 3, it could be possible to obtain some properties of
(moderate) growth on the $F_i$'s.\\
{\bf 7.} We can assume that equation ${\sf L}{ \sf
  de}_{i}$ is ``totally reduced'' i.e. that another application
  of one step of Abel's method gives a null equation. In this case, it is
  interesting to study the quotient $  
 {\underline{{\cal S}\scriptstyle{\stackrel{ \cal O}{\omega_i}}}( {\sf
    L}{ \sf de}_{ i})/     \left[ {\underline{{\cal
     S}\scriptstyle{ \stackrel{\cal O}{{}_\omega}}}} \right]_i }$. We
  conjecture that it is trivial. Combined with remark 6, this
  could give supplementary informations about the nature of
  equation $ ({\sf L}{\sf de}_i)$.

\subsection{{\bf Two Methods to solve {\sf Afe} with real rational
    inner functions}}
The proof of the preceding theorem contains some useful tools to construct
two ``methods'' of solving {\sf Afe} of  the type $ \eq{E}
$.\\
The first tool is essentially based on Abel's method. It is well
formalized and appears 
very general: its only defect is to be computational.\\
The second only consists in a remark and is not well established as a
general method. Meanwhile, this remark allows us in part 3 to solve 
the two $ {\sf Afe} $  $ \eq{R} $ and $ ( {\cal SK})$, associated
respectively to equation $(L_2)$ and Spence-Kummer equation
$(SK)$. \\
It is based on the idea that certain
solutions of $ \eq{E} $ are determined by their monodromy.
\subsubsection{{\bf Abel's method of resolution of AFE with rational
  inner functions}}
Let us assume that $ (U_1,...,U_N) \in \rR (x,y)^N $ is such that there exists a non-constant holomorphic genuine solution $ {\bf F}=(F_1,..,F_N) \in 
 \underline{{\cal
      S}{ \stackrel{\cal O}{{}_\omega}}}   ( \cal E) $.
Then let $ \left [ {\underline{{\cal
     S}\scriptstyle{ \stackrel{\cal O}{\omega}}}} \right]_i $ be the subspace of holomorphic germs at $\omega_i$ spanned by the $i$-th components of
 solutions $  {\bf F} \in \underline{{\cal
      S}\scriptstyle{ \stackrel{\cal O}{\omega}}}   ( \cal E) $.\\
We can choose $\omega \not \in \Sigma_{\sf R}$ such that 
for each
 $i\in \{1,..,N\}$, there is  a non-trivial linear differential equation $ ({\sf L}{ \sf de}_{ i})$ having
algebraic coefficients which are well defined at $\omega$. This equation is such
that every component $F_i$ of~${\bf F} \in \underline{{\cal
      S}\scriptstyle{ \stackrel{\cal O}{\omega}}}   ( \cal E) $ satisfies $ ({\sf L}{ \sf de}_{ i})$. We note $
  \underline{{\cal S}\scriptstyle{\stackrel{ \cal O}{\omega_i}}}( {\sf
    L}{ \sf de}_{ i}) \supset  \left[ {\underline{{\cal
     S}\scriptstyle{ \stackrel{\cal O}{\omega}}}} \right]_i $ the linear space of
 the holomorphic germs at $\omega_i$ which are solutions of this equation. Let  $
 \{ {\sf G}_i ^{\nu} \, | \, \nu =1,..,\nu_i \, \} $ be a
 basis of this space.\\
Then we have
$$    \underline{{\cal
      S}{ \stackrel{\cal O}{{}_\omega}}}   ( {\cal E})=
  \left\{ ( \sum_{\nu=1}^{{\nu}_1} a_1^{\nu} {\sf G}_1^{\nu},...,  
\sum_{\nu=1}^{\nu_N} a_N^{\nu} {\sf G}_N^\nu) \in 
\bigoplus_{i=1}^N  \underline{{\cal S}\scriptstyle{\stackrel{ \cal
      O}{\omega_i}}}( {\sf L}{ \sf de}_{\, i}) \, | \, \sum_{i=1}^N
\sum_{j=1}^{\nu_i} a_i ^{\nu} {\sf G}_i ^{\nu} ( U_i)=0 \, \right\} $$ 
so, in a certain way, the explicit resolution of $ \eq{E} $ at $\omega $
amounts to some linear algebra in a finite dimensional space.\\

It is easy to prove that, in the standard coordinates system $(x,y)$ on
$\rC^2$, the derivations ${\partial}_p^{kl}$ (where $ p=l,k$) are elements of 
 $ \rC(x,y) \partial_x+\rC(x,y) \partial_y $ .\\
 Then the  coefficients ${\cal A}_{ij} $ of any {\sf Adfe} obtained
 through the application of several steps of Abel's method to $ \eq{E} $
 belong to $ \rC(x,y)$, therefore the process to obtain $ ({\sf
   L}{\sf de}_i)$ from the {\sf Afe} $\eq{E} $ can be performed
 within $ \rC(x,y)$. This fact allows us to easily implement an
 algorithm on a computer algebra system which constructs $ {\sf
   L}{\sf de}_1$ from $ (U_1,U_2,..,U_N)$  . \\
The author has used this method to solve the equation $ ({\cal
  E}_{\sf c})$ of part 3.5, and it seems possible to apply it to all
the equations of part 3 .

\subsubsection{Method of monodromy ``a priori''}

 Contrarily to the preceding method, the one described here doesn't
seem to be valid in the general case, but its interest lies in the fact that it works
for at least three ${\sf Afe}$ associated to classical functional
equations of polylogarithms $\l{k} $ with $ k \leq 3$ (see part 3).\\
It is a ``method'' to find  solutions with logarithmic growth of an
{\sf Afe} when the $U_i$'s verify a certain condition called
``condition $ ({ C})\,  $'', which is defined below. Roughly speaking, it is based on the
fact that solutions with logarithmic growth are determined by their
monodromy, which can be determined ``a priori'' when the solutions of
some sub-equations of $ \eq{E} $ are known.\\

We now define ``condition $({ C})\, $'' : 
 \begin{definition}
  {\em The set of rational functions $ \{ \, U_i \, \} $ verifies
  `` condition $ ({\cal C}) $~'' if for all $ i \in \{ 1,..,N \} $
  there exists $ l(i) \neq i $ such that $(U_i,U_{l(i)} ) $ is a global
  system of coordinates on $ X:= \projc ^2 \setminus \Sigma $ .}
 \end{definition} 
In the following pages, we will assume this strong condition verified.\\

Let  $ {\bf F}=(F_1,..,F_N) $ be a genuine solution of $ \eq{E} $ at a
generic point $ \omega \in \rR^2 \setminus \Sigma$.\\
Let be $ i \in \{1,..,N\}$. There exists an integer $\sm_i  $ and a finite
number of distinct points $ a_k^i\in \projc ^1$ ($1\leq k \leq \sm_i$
) such that $ X_i=U_i( X)= \projc^1 \setminus \{
a_{\nu}^i \, \}_{ \nu \leq \sm_i \, }$. Theorem 1 implies that
every germ $F_i$ at $ \omega_i$ can be analytically
extended along any path in $ X_i$.\\

Let be $ \Lambda_i= \{ \gamma
^{\lambda}_i \} _{\lambda \leq \sm_i} $ a minimal family of loops of
basepoint $\omega_i$ in $
X_i $, such that their homotopy classes and
their inverse span $ \Pi _1(X_i,\omega _i) $ (a suitable choice is to take for $\gamma_i^{\lambda} $ a loop in $X_i$, of index $1$ with
 respect to $a_j ^{i}$ if $j=\lambda$, and of index 
$0$ otherwise) .\\

Now we fix $i$ and we note $l$ for $l(i)$. Condition $ ( C)$
implies that we can find a loop $\overline{\gamma}
^{\lambda} _{i} $ of basepoint $\omega$ in $X$ such that 
$ [U_l \circ \overline{\gamma}
^{\lambda} _{i}]=[1] $ in $ \prod_1(X_l,\omega _l)$ and $ 
[U_i \circ \overline{\gamma}
^{\lambda} _{i}]=[\gamma ^{\lambda}_i] $ in $ \prod_1(X_i,\omega
_i)$ .\\

Because we have  $ F_1(U_1)+F_2(U_2)+....+F_N(U_N)=0$ in  a
neighbourhood of $\omega$, then by analytic continuation along 
$\overline{\gamma}^{\lambda} _{i} $ we get a new functional relation 
in $ \germhol{{\omega}} $
: 
$$ F_1^{[U_1 \circ\overline{\gamma}
^{\lambda} _{i}]}  (U_1)+...+ F_i ^{[{\gamma}
^{\lambda} _{i}]}   (U_i)+...  + F_l ^{[1]}(U_l)+...
.+F_N ^   {[U_N \circ\overline{\gamma}
^{\lambda} _{i}]}          (U_N)=0 $$
which can be summarized by $ {\bf F}^{[\overline{\gamma}           ^{\lambda}_{i}]} \in  
 \underline{{\cal
      S}{ \stackrel{\cal O}{{}_\omega}}}   ( \cal E) $ where 
$ {\bf F}^{[\overline{\gamma} ^{\lambda}_{i,j}]}:=( F_k ^{[U_k \circ \gamma
  ^{\lambda}_{i}]} )_{k=1..N} $ .\\

By taking the difference between the above equations, we get a new one 
\begin{align} (  F_1^{[U_1 \circ\overline{\gamma}
^{\lambda} _{i}]}  (U_1) -F_1(U_1) ) &  +...  +(  F_i ^{[{\gamma}
^{\lambda} _{i}]}   (U_i)-F_i(U_i) ) +..      \nonumber \\ 
..  &  + ( F_l
^{[1]}(U_l)-F_l(U_l ) +...
.+( F_N ^   {[U_N \circ\overline{\gamma}
^{\lambda} _{i}]}          (U_N)-F_N(U_N) ) =0  \nonumber
\end{align}
Now the germ $ F_l ^{[1]} -F_l$ is null in $ \germhol{ \omega_l} $,
therefore $ {\bf F}-{\bf F}^{[\overline{\gamma} ^{\lambda}_{i}]} $ is not a
genuine solution of $\eq{E} $ any more.\\

Let be $ K_ i ^{\lambda} =\left\{ k \, |   \,[U_k\circ
\overline{\gamma}^{\lambda} _{i}] \neq  [1] \mbox{ in } \Pi_1(X_k,\omega_k)
\, \right\} $. We have $ K_ i ^{\lambda}   \varsubsetneq \{1,..,N\}$.\\
Let us assume that we know a basis $ \{ {\bf
B}_{i,\kappa}^{\lambda} \, | \, \kappa \in \Delta_i ^{\lambda} \}  $
of  $  \underline{{\cal S}}_{K _ i ^{\lambda}   } $ , with 
$ {\bf B}_{i,\kappa}^{\lambda} =({\bf b}_{i,\kappa}^{\lambda , 1}, ..,
{\bf b}_{i,\kappa}^{\lambda , N}) $.

Then we get a relation
$$  {\bf F}-{\bf F}^{[\overline{\gamma} ^{\lambda}_{i}]} = \sum_{ \sigma \in
  \Delta_{i}^{\lambda}}  {\beta}^{\lambda }_{i, \sigma} \,   {\bf B}_{i ,
\sigma} ^{\lambda }
\qquad \mbox{ with  } \:
{\beta}^{\lambda }_{i, \sigma} \in \rC \qquad $$
from which we get the following relations for all $\lambda \in \Lambda_i $ : 
\begin{align}
 \qquad {\cal M}_{[\gamma ^{\lambda}_i ]} \, F_i & = F_i +  \sum_{ \sigma \in
  \Delta_i ^{\lambda} }  {\beta}^{\lambda }_{i, \sigma}     \,  {\bf b}_{i ,
\sigma} ^{\lambda  , i } 
  \qquad  && (\star)_{i}^{\lambda} \nonumber \\
\qquad {\cal M}_{[U_s \circ \overline{\gamma} ^{\lambda}_{i} ]} \, F_s & =
  F_s +  \sum_{ \sigma \in \Delta_i ^{\lambda} }
      {\beta}^{\lambda }_{i, \sigma}   \,  {\bf b}_{i ,
\sigma} ^{\lambda  , s } 
  \quad 
  \mbox{for}  \: s \in   K_ i ^{\lambda} \: \mbox{    and    }\:  s
  \neq i  && (\star  \star )_{i,s}^{\lambda} \nonumber 
\end{align}

If $Y$ is a complex manifold, knowing the monodromy of $ G\in
\tildhol{Y} $ means knowing a representation
\begin{align}  {\Pi} _1( Y,y) & 
\longrightarrow   \mbox{ End}_{ \scriptstyle{\rC}} \bigl( \, 
\underline{Dy} \, \bigr)  \nonumber \\
[\gamma] & \longrightarrow  \quad T_ {[\gamma]} \, : g \rightarrow g^{[\gamma]}
\nonumber 
\end{align}
for at least one $  y \in Y $, where $ \underline{Dy} $ denotes the
linear space of the 
determinations of $G$ at $y$ .\\

Because we have chosen the family  $ \{ \, [ \gamma _i ^{\lambda } ] \, ,
\, [ \gamma _i ^{\lambda } ]^{-1} \, \}$ such that it spans $ \prod_1
( X_i , \omega_i) $ , the relations $ (\star)_{i}^{\lambda}$ give us `` a priori '' the
monodromy of each of the components $F_i$ in function of the
components $  {\bf b}_{i ,\sigma} ^{\lambda  , s } $ of the
subsolutions   $ {\bf B}_{i ,
\sigma} ^{\lambda }  $ of $ \eq{E} $ (for $i=1,..,N$ and $ \lambda
\leq {\sf m}_i$).
\begin{prop}
Under condition (C), the monodromy of each of the components $F_i$ of a genuine solution of $
 \eq{E}  $ can be expressed in terms of the
 components of some subsolutions of $ \eq{E} $.
\end{prop} 
This transforms our point of view on equation $
\eq{E} $ : although considering it in a functional form we will now 
 see relations $ (\star) $  as  {\em `` monodromy equations''}
for the components of the solutions, and  relations $ (\star \star) $ as  
{\em `` compatibility relations '' } between those equations of 
monodromy.\\ 

We now want to find some genuine solution of $ \eq{E} $ by 
``solving'' the monodromy equations~$(\star) $.  \\ 

Let us assume that there exists a genuine solution ${\bf F}=(F_1,..,F_N)$ of
$\eq{E} $ at $\omega$.\\
 From the preceding lines, it comes that there
exist complex constants $ {\beta}^{\lambda }_{i, \sigma} ( {\bf F})$
satisfying both relations $(\star)$ and $ (\star \star)$. \\

Let be $\widetilde{\bf H}=(\widetilde{H}_1,..,\widetilde{H}_N) \in
\prod_i 
{{\widetilde{\cal
      O}}}( X)$ such that each $ \widetilde{H_i}$ has a determination
      $ H_i$ at $ \omega_i$ satisfying the equations $ (\star)_i
      ^{\lambda} $. Then the germ $H_i-F_i$ can
      be extended analytically to $X_i$ without ramifications. This implies that the germ $ {\cal H}=\sum (H_i-F_i) \circ U_i $ at
      $\omega$ develops into a global holomorphic function\,:  $ {\cal
      H}\in {\cal O}(\projc^2
      \setminus \Sigma)$.\\

Now let us suppose that we can choose $ \widetilde{H_i}$ with logarithmic
growth. Then $ {\cal H} $ is a global holomorphic function on $
\projc^2\setminus  \Sigma$ with logarithmic growth at infinity. 
By a Liouville type theorem, this implies that ${\cal H}$ is constant. 
Then ${\bf
  H}=(H_1,..,H_N) \in  \underline{{\cal
      S}{ \stackrel{\cal O}{{}_\omega}}}   ( \cal E)$. \\

We note $\underline{{\cal
      S}\scriptstyle{ \stackrel{{\cal O}}{{}_\omega}}}   ( {\cal
      E})^{\scriptstyle{log}}$ the subspace of the solutions of $ \eq{E} $
      with logarithmic growth.\\ 
Then the problem of finding genuine solutions in  $\underline{{\cal
      S}\scriptstyle{ \stackrel{{\cal O}}{{}_\omega}}}   ( {\cal
      E})^{\scriptstyle{log}}$ amounts to solving the equations~$(\star)$ in
the space $ \prod_i {\underline{ {\cal O}\scriptstyle{\stackrel{ log}{
    \omega_i}}}}$ $ (X_i)$. One of the conceptual interests of this
      is that the problem is now reduced into a linear form.\\

Let be ${\bf F} \in \underline{{\cal
      S}\scriptstyle{ \stackrel{{\cal O}}{{}_\omega}}}   ( {\cal
      E})^{\scriptstyle{log}}$. Then the subsolutions of the form 
$  {\bf F}-{\bf F}^{ [ \overline{\gamma}_i^{\lambda}]} $ which appear in the
      preceding discussion are now elements of $ \underline{{\cal
      S}\scriptstyle{ \stackrel{{\cal O}}{{}_\omega}}}   ( {\cal
      E}_{K_i ^{\lambda}})^{\scriptstyle{log}}$. Under suitable
      conditions on the $U_i$'s, the $ \{ U_j \}_{j\in
      {K_i ^{\lambda}}} $ verify condition $ ({ C}) $ again.
In this case it could be possible to inductively determine 
   the solutions with logarithmic growth of equation $ \eq{E} $.\\ 

 {\bf remarks} {\bf 1.} Most components of most of the solutions of known
$\eq{E} $-form equations are constructed from iterated
integrals (see part 3). Then it will appear interesting and 
useful to work in the subspace $ \prod_i \underline{{\cal I}\scriptstyle{ \omega_i}}
\varsubsetneq \prod_i {\underline{ {\cal O}\scriptstyle{\stackrel{ log}{
    \omega_i}}}}$ $ (X_i)$ where $ \underline{{\cal I}\scriptstyle{ \omega_i}} $
denotes the space of the determinations of the elements of $
{\cal I}_{ \{X_i\}}$ at $\omega_i$. \\
{\bf 2.} But not all the components of the solutions with logarithmic
growth are constructed from iterated integrals: for instance the
function  $ \mbox{Arctan}( \sqrt{\bullet})$ is a
component of a solution of the ${\sf Afe}$ with real rational inner
functions $(\cal SK)$ considered in 3.4. This function cannot be
expressed from iterated integrals although it is ramified with
logarithmic growth on $\projc^1 \setminus \{ 0,1,\infty\} $.

\section{Examples of explicit resolution of abelian functional equations
  with real rational inner functions}
In this part we apply the method sketched above to the resolution of
some functional equations: to begin with, we solve some very classical
equations which have already been 
treated by Abel using his own method in \cite{abel}. Here we use some
monodromy arguments to
solve them. We finish
with the ``generalized  Spence-Kummer equation of the trilogarithm''  
and with another one which will be interpreted  in the
framework of web geometry in part~4.1

\subsection{Cauchy equation revisited}
Here we want to solve again the ``generalized Cauchy equation '' $ \eq{C} $ in 3 unknowns 
$$   F_1(x)+F_2(y)+F_3(\frac{x}{y})=0  \qquad \eq{C} $$
by using monodromy arguments: we are interested in solutions the

We note $U_1(x,y)=x, \; U_2(x,y)=y, \;  U_3(x,y):=\frac{x}{y}
$,  and $ {\cal W}_{\cal C} $ the web given by the three foliations, the leaves of
which are respectively the level curves of $U_1, U_2$ and $ U_3 $ .
Its singular locus is $\Sigma_{\cal C} := \{ (z,\zeta)\in  \rC ^2 \:
| \;  \;  z \; \zeta=0  \; \} $ . \\

An easy computation gives us that $ U_i( \rC ^2 \setminus \Sigma_{\cal
  C} ) = \rC ^{\star} $ for $ i=1,2 \mbox{  and  } 3 $ . 
We deduce that if $ (F_1,F_2,F_3) \in \underline{{\cal
      S}{ \stackrel{\cal O}{{}_\omega}}}   ( \cal C) $ where
  $\omega=(1,1)\not \in \Sigma_{\cal C} $, then $ F_i \in
   \widetilde{{\cal O}}({ \rC ^{\star} })  $ for $i = 1,2, \mbox{  and
     } 3 $ .\\

In this case, condition $(C)$ is verified. We are looking for the
solutions of $ \eq{C} $  the components of which are elements of the
space $ {\bf {\cal I}}_{\{0\}} $ of the iterated
integrals on $ \projc ^1
\setminus\{0,\infty\} $ relative to the rational $1$-form $
\omega_0:={dz} / {z} $ : we have  $ {\bf {\cal I}}_{\{0\}} :=
\mbox{Vect}_{{}_{\scriptstyle{\rC}}} \langle\;  \{ \;  \log^k(\bullet)
\}_{\scriptstyle{k\in \rN}} \rangle $~. \\

Let $\gamma_0$ be a loop with $1$ as its base point turning around $0$ in
the direct sense\,: its homotopy class $[ \gamma ] $ is a generator of $
\Pi_1 ( \rC ^{\star} , 1 ) \simeq \rZ $ .\\
There are two loops in $ \projc ^2 \setminus \Sigma_{\cal C} $,
noted $\gamma ^1$ and $\gamma ^2 $, such that we have 
$ U_i\circ \gamma ^j = \gamma $ if $i=j$, and $ U_i\circ \gamma ^j$ is the constant path otherwise.
By analytic continuation of $ \eq{C} $  along $\gamma_1 $ we get a new
functional equation. By taking the difference between these
two equations, it comes that in a neighbourhood of $\omega$, we have
$$ ( F_1^{[\gamma]}(x)-F_1(x)) +( F_3^{[\gamma]}(\frac{x}{y})-F_3(\frac{x}{y}))
=0 $$  
Both this equation and the one given by analytic continuation along $
\gamma ^2 $ imply that there exists a constant $a\in\rC $ such that 
$$   {\cal M}_0 \, F_1 =  F_1 +a \quad , \quad 
 {\cal M}_0 \, F_2 = F_2 +a \quad , \quad
 {\cal M}_0 \, F_3 =  F_3 -a $$
Considering these relations as equations of monodromy in the algebra $
{\bf {\cal I}}_{\{0\}} $, we get only one possible solution (modulo
the constants)
$$ {\bf L}:= a \;( \int^{\,  \bullet} \omega_0 , \int^{\, \bullet} \omega_0,-
\int^{ \, \bullet} \omega_0  )=a \; ( \llog(\bullet) ,\llog(\bullet),- \llog(\bullet)) $$
Using Bol's bound of proposition 4, we obtain that, for all $\widehat{\omega} \not \in
\Sigma_{\cal C}$, modulo the constant solutions,
$\underline{{\cal S}{ \stackrel{\cal O}{{}_\omega}}}   ( \cal C) $ is
spanned by any analytic
continuation of $ {\bf L} $
from $\omega$ until~$\widehat{\omega}$~in $\projc^1 \setminus \Sigma_{\cal C} $.  

\subsection{Arctangent  equation revisited(see \cite{abel})}

It is well known that the arctangent function $ {\bf A}rc(\bullet) :=
\int_0 ^{ \, \bullet } \frac{dx}{1+x^2} $ satisfies the functional
equation 
$$  \qquad  \qquad  \qquad {\bf A}rc ( x ) + {\bf A}rc ( y) = {\bf A}rc (\frac{x+y}{1-xy} )
\qquad  \qquad  \qquad  ({\cal A}rc) $$
on the two real sets, $ \{ \, xy<1 \, \} $ and $ \{ \, xy>1 \, \} $.\\
We consider a generalized version of $ ({ Arc})$ (with 
$V_1(x,y)=\frac{x+y}{1-xy}$ )
$$ \qquad \qquad \qquad G_1(U_1)+G_2(U_2)+G_3( V_1 )=0  \qquad \qquad
\quad \qquad ({\cal A} rc) $$
Then the singular locus of
the web associated to $ ({\cal A} rc) $  is 
 $$ \Sigma_{\scriptstyle{
    {\cal A}rc}}:= \{  (z,\zeta)\in  \rC ^2 \:
| \;  \;  (1-z \; \zeta) \, (1+z^2)\, (1+ \zeta ^2)=0  \; \} $$ 
Let be $\omega:=(0,0) \not\in  \Sigma_{\scriptstyle{\cal A}rc }$ 
. We want to determine $\underline{{\cal
      S}\scriptstyle{ \stackrel{\cal M}{\omega}}}   ( {\cal A}rc)  $ .\\
If ${\bf A}:=(A_1,A_2,A_3)\in \underline{{\cal
      S}\scriptstyle{ \stackrel{\cal M}{{}\omega}}}   ( {\cal A}rc) $, then, in the
  same way than in 3.1, we get that the $A_i$'s
  are global analytic functions ramified in $ +i$ and $-i$ : 
$ A_j \in \underline{ \widetilde{\cal O}\scriptscriptstyle{\omega}} ( \rC \setminus \{ \pm i \})  $ .\\
We are looking for solutions, the components of which are elements of
the algebra ${\bf {\cal I}}_{ \{ \pm i \} } $.\\

By using the method of monodromy ``a priori'' 
we obtain the following relations of monodromy for the $A_i$'s:
\begin{align}
{\cal M}_i \, A_1 & =A_1 + a  && {\cal M}_{-i} \, A_1  =A_1 -a
\nonumber \\ 
 {\cal M}_i \, A_2 & =A_2 + a   && {\cal M}_{-i} \, A_2  =A_2 -a \quad 
 \quad   ({\cal M } o)  \nonumber                    \\ 
{\cal M}_i \, A_3 & =A_3 - a   && {\cal M}_{-i} \, A_3  =A_3 +a
\nonumber 
\end{align}
where $a \in \rC $ is a constant. \\
The relations $({\cal M}o)$ considered as equations in ${\bf {\cal
    I}}_{ \{ \pm i \} }^3 $ admit a single possible solution (modulo
the constants):
$$ A_1( \bullet ) = A_2( \bullet ) = -A_3( \bullet ) =a\, \int_0^{\,
  \bullet } \omega_i - a \, \int_0^{\, \bullet } \omega_{-i}= 2ia\,
\int_0^{ \, \bullet }  \frac{dz}{1+z^2} $$ 
For dimensional reasons, we obtain that, for all $\widehat{\omega} \not \in
\Sigma_{\cal C}$, modulo the constant solutions,
$\underline{{\cal S}\scriptstyle{ \stackrel{\cal O}{\omega}}}   ({\cal A}rc ) $ is
spanned by any analytic
continuation in $\projc^1 \setminus \Sigma_{\cal C} $ of $ ( {\bf
  A}rc , \, {\bf A}rc , \, {\bf A}rc \, )  $
from $\omega$ until~$\widehat{\omega}$~.

\subsection{Roger's dilogarithm equation revisited (see \cite{blabol},  \cite{roger})}

 In \cite{roger}, L. Rogers established a ``clean version'' of the equation
 $ (L_2)$ verified by the Rogers dilogarithm ${\bf d}$,  for $0<x<y<1$\,:
$$ \qquad  \qquad 
 {\bf d}(x)- {\bf d}(y)- {\bf d}(\frac{x}{y})- {\bf d}(\frac{1-y}{1-x})+
 {\bf d}(\frac{y(1-x)}
{x(1-y)})=0 \qquad  \qquad   \qquad \ ({ R})$$
(here we have taken  ${\bf d}( \bullet) := \l{2} (
\bullet) +\frac{1}{2} \llog( \bullet ) \llog(1-\bullet ) -\frac{ \pi
  ^2}{6} $ : it is a normalized version of the original Rogers dilogarithm
(by addition of $-\pi ^2/6$) in order to have $0$ for the rhs of ($R$)). \\

We consider the more general equation in 5 unknowns\,:
$$  \qquad \qquad  D_1(x)+ D_2(y)+D_3(\frac{x}{y})+ D_4(\frac{1-y}{1-x})+
 D_5(\frac{y(1-x)}{x(1-y)})=0 \ \qquad   \qquad  \ ({\cal R})$$
We note ${\cal W}_{ \scriptstyle{ {\cal R}}}$ the singular web associated
to the inner functions $U_1,U_2,...,U_4,U_5 $ of  $ ({\cal R})$  , where
$U_4(x,y):=\frac{1-y}{1-x}$ and $U_5(x,y):=\frac{y(1-x)}{x(1-y)}$ .\\
After computation we get that its singular locus is 
$$ \Sigma_{ \scriptstyle{ {\cal R}}}:= \{ (z,\zeta)\in  \rC ^2 \:
| \;  \;  z \; \zeta \, (1-z) \, (1-\zeta) (z-\zeta)=0  \; \} $$ 
We choose  $\omega:=(\frac{1}{3}, \frac{1}{2}) \in \rR ^2 \setminus 
 \Sigma_{ \scriptstyle{ {\cal R}}}$.\\
In \cite{bol}, G. Bol found an equivalent of a basis of this space: in the
framework of web geometry (see part 4.1. below), he
determines a basis of the space of abelian relations of ${\cal W}_{
  \scriptstyle{ {\cal R}}}$. 
We want to rediscover Bol's results by application of our two
``methods'' described in part 2.5. \\

{\bf {\sf A)} Resolution of $ \eq{R} $ by the method of ``monodromy a
  priori''}\\

By an easy computation we find that $ U_i( \rC ^2 \setminus  \Sigma_{
  \scriptstyle{ {\cal R}}} ) = \rC \setminus \{ 0,1 \} $ . So, if 
${\bf D}=(D_1,..,D_5)\in \underline{{\cal S}{ \stackrel{\cal
      O}{{}_\omega}}}   ( {\cal R})$, then $ D_i\in  \underline{ {\cal
  O}\scriptstyle{\omega}}({ \rC
  \setminus \{ 0,1 \} }) $ for $ i=1,..,5$ .\\
In this case, equation $ \eq{R} $ can be solved by the method of
  monodromy a priori. We want to determine the solutions of $ \eq{R} $
  the components of which are iterated integrals elements of $ {\cal
  I}_{ \{0,1\} }$.\\
We begin to search the $3$-solutions of this type: we want to determine 
$ 
F^3 \underline{  {\cal S} \scriptstyle{ \stackrel{{\cal
  I} }{\omega}   }}   ( {\cal R})$.\\
Our method of ``monodromy a priori'' works very well without
difficulties and too many computations.
It gives us the following $5$ non-constant  independent elements of $ F^3
\underline{{\cal S}{ \stackrel{\cal I}{{}_\omega}}}   ( {\cal R})$:

\begin {align}
 &{\bf \Delta_1}:= \biggl(    \:  {\bf L}_{x_0}, -{\bf L}_{x_0},-
  \: {\bf L}_{x_0}\: , \: 0 \: ,0   \biggr) 
&&  {\bf \Delta_2}:= \biggl( \: 0 \: , \: 0 \: , \: {\bf L}_{x_0} , {\bf L}_{x_0} ,   - {\bf L}_{x_0} \:  \biggr) \nonumber  \\ 
 &{\bf \Delta_3}:= \biggl(    \: {\bf L}_{x_1}, -{\bf L}_{x_1},0,-{\bf
  L}_{x_0} \:,0 \:   \biggr)    &&   {\bf \Delta_4}:=\biggl( \:  {\bf L}_{x_1}\: ,0 \: , 
- {\bf L}_{x_1} \: ,\: 0 \: ,\: {\bf L}_{x_1} \: \biggr)
     \nonumber \\
& {\bf \Delta_5}:=\biggl( \: {\bf L}_{x_1+x_0},0\: ,
 -{\bf L}_{x_1+x_0}\:  ,   {\bf L}_{x_1} \: , 0  \biggr)    && \nonumber  
\end{align}

Now we have to try to determine the last non-constant solution of
$\eq{R} $ if there is one. We will use our method again\,:
we want to detail the computation to be well understood.

Let us consider the loop $\gamma: \, [0,1] \in \sigma \rightarrow (\exp(2i\pi
  \sigma)/3,1/2) \in  \rC ^2 \setminus  \Sigma_{ \scriptstyle{ {\cal R}}}  $.
The computations give 
\begin{align}
[U_1\circ \gamma ] & = [ c_0^1]    &  [U_2\circ \gamma] & = [1]    &
[U_3\circ \gamma]& =[ c_0^3]  \nonumber  \\
[U_4\circ \gamma]& =[1]   & [U_5\circ \gamma ] & =[c_0^5]  &  \nonumber  
\end{align}
where these equalities are (respectively) in $ \Pi_1( \rC \setminus \{ 0,1 \} ,
\omega_i) $ , for $ i=1,..,5 $ .\\
So we have a new functional equation
$$   ( D_1^{[c_0^1]}(U_1)-D_1(U_1))+  ( D_3^{[c_0^3]}(U_3)-D_3(U_3))+
 ( D_5^{[c_0^5]}(U_5)-D_5(U_5))=0$$
which corresponds to an element of  $ F^3 \underline{{\cal S}{ 
\stackrel{\cal O}{{}_\omega}}}   ( {\cal R})$ .
But we explicitly know this space, and in this case we obtain
the following relations of monodromy for the components of any
solution ${\bf D}$: 
\begin{align}
{\cal M}_0 \, D_1 & = D_1 +a \,{\bf L}_{x_1} +a_1 \nonumber \\
{\cal M}_0 \, D_3& = D_3 -a \,{\bf L}_{x_1} +a_2 \nonumber \\
{\cal M}_0 \, D_5& = D_5 +a \,{\bf L}_{x_1} -(a_1+a_2) \nonumber 
\end{align}
where $ a,a_1,a_2,  a_3 $ are complex constants.\\

Now considering the path $ \sigma \in [0,1]  \rightarrow 
( \frac{1}{3} , 1- \frac{1}{2} \exp (2i\pi\sigma )\; ) $ in 
$\rC ^2 \setminus \Sigma_{ \scriptstyle{ {\cal R}}} $ 
we get by the same way 
\begin{align}
{\cal M}_1 \, D_2 & = D_2 +{a}' \, {\bf L}_{x_0}+a_1' \nonumber \\
{\cal M}_0 \, D_4 & = D_4 +{a}' \, {\bf L}_{x_1}+a_2' \nonumber \\
{\cal M}_0 \, D_5  & = D_5 -{a'} \,{\bf L}_{x_1} -(a_1'+a_2') \nonumber 
\end{align}
From these relations it comes that $ a=-a'$ and  $ a_1+a_2= a_1'+a_2'$ .\\
We can continue this type of computation and finally we get that the
monodromy ``a priori'' of the components of holomorphic solutions of $
\eq{R} $ are ($a_i,b_j $ being complex constants satisfying certain
linear relations)
\begin{align} {\cal M}_0 D_j &  = D_j- \epsilon_j \, a \, 
{\bf L}_{x_1} +  a_j \, \nonumber \\
 {\cal M}_1 D_j&  = D_j+ \epsilon_j \, a \, 
 {\bf L}_{x_0}+ b_j  \nonumber
\end{align}
with $\epsilon_j=1 $ for $ j=1,5$ and $-1$ otherwise.
It can be proved that the $a_i$'s and $b_j$'s are such that there
exists a linear combination $ {\bf H}=\sum \alpha_i {\bf \Delta}_i$ such that
the monodromy of the components of $ {\bf D'}=(D_j'):={\bf D }+{\bf H}$
verifies
\begin{align} {\cal M}_0 D_j' &  = D_j'+ \epsilon_j \, a
   {\bf L}_{x_1}    \nonumber \\
 {\cal M}_1 D_j' &  = D_j'+ \epsilon_j \, a
  {\bf L}_{x_0}  \qquad \qquad \qquad  \nonumber
\end{align}
Now it is not difficult to prove that the function $ {\bf f}=a {\bf
  d} \in {\cal I}_{ \{0,1\} }$  satisfies the following 
monodromy equations 
$$   
{\cal M}_0 {\bf f}   ={\bf  f}-a \, {\bf L}_{x_1} \: , 
 \quad {\cal M}_1 {\bf f}   = {\bf f}+a \, {\bf L}_{x_0} 
 $$
We deduce that $  {\bf \Delta}_6:= ( {\bf d}+c, -{\bf d} ,-{\bf d}, -{\bf d},
{\bf d} ) \in   \underline{{\cal S}{ 
\stackrel{\cal O}{{}_\omega}}}   ( {\cal R})$ , where $ c $ is
a constant (in fact $ c =0$). \\
We can easily construct a basis $ \{ {\bf \Delta }_i \, | \, i=-3, ..,0 \,
\} $ of the constant solutions of $  \eq{R} $ .
It is not difficult to prove that the $10$ elements $ {\bf \Delta}_j $
described above are linearly independent.\\
On the other hand  we have 
$$ 10 =     \mbox{ dim} _{ \scriptstyle{ \rC}}  \, \left< \, \{ \, {\bf \Delta}_j \,
\} \, \right> \leq 
\mbox{ dim} _{ \scriptstyle{ \rC}} \, \underline{{\cal S}{ 
\stackrel{\cal O}{{}_\omega}}}   ( {\cal R})   \leq {5(5-1)}/{2}=
10 $$
where the last inequality is given by proposition 4. Then we deduce that 
$$  \underline{{\cal S}{ 
\stackrel{\cal O}{{}_\omega}}}   ( {\cal R})  = \left< \{ {\bf \Delta}_j \,
| \, {j=-3,-2,..,6}  \, \} \, \right> $$ 
This solves $ \eq{R} $ at $ \omega$ in the holomorphic class.
We get the local holomorphic solutions around $ \omega'\not \in
\Sigma_{\cal R} $ by analytic continuation of the $ {\bf \Delta}_j $'s
along any path joining $ \omega $ to $ \omega' $ in~$X$~. \\                                                                                
    
{\bf {\sf B)} Resolution of $ \eq{R} $ by Abel's method }\\
A simple application of Abel's method implies that on the
whole $\Omega$, the first component of every solution of $ {\cal R}$ must verify the following linear differential
equation 

$$     \frac{d^4 g }{dv^4}
+\frac{4(2v^3-3v^2+v)}{ v^2(1-v)^2}  \,  \frac{d^3 g }{dv^3} +
\frac{2(1-7v+7v^2)}{ v^2(1-v)^2}     \,  \frac{d^2 g }{dv^2}      
+\frac{2(2v-1)}{ v^2(1-v)^2} \, \frac{d g }{dv}       =0     
$$
By integrating this equation, which can be done without great
difficulty with a
computer system, we find that it admits as 
 general solutions the functions of the form $ c_1\, {\bf
  d}+ c_2 \, {\bf L}_{x_0} + c_3 \, {\bf L}_{x_1}   +c_4 $, 
which is the form that any first component of any
solution of $\eq{R} $ can have. \\

{\bf remark:} but even without integrating, the last equation 
gives us some informations:
it admits three singular points, $0,1$ and $\infty$. One can easily prove
 that they are regular points. By a classical theorem of
the theory of linear differential equations with mereomorphic
coefficients, it comes that any solution of $ ( {\cal R}_{\sf de}^4)
$ a priori has moderate
growth near~$0,1$~and~$\infty$. Another remark is that the differential
operator associated to this equation can be factorized into a product of
differential operators of first order.

\subsection{Spence-Kummer equation of the trilogarithm visited(see
  \cite{lewin})}
To the Spence-Kummer equation $(SK)$ satisfied by $\l{3} $ ( with $0<x<y<1$ ) 
we can associate the following abelian functional equation 
$$  \qquad \quad \qquad  \qquad 
F_1(U_1)+F_2(U_2)+F_3(U_3)+...+F_9(U_9)=0   \qquad \qquad
\qquad  ({\cal SK})
$$
where the $U_i$'s are the rational inner functions which appear in $({ S
  K})$: $U_1,U_2,...,U_5$ have been defined above and we note 
\begin{align}
 U_6(x,y) & =xy   && U_7(x,y) 
 = \frac{x(1-y)}{x-1}  \qquad \qquad  \nonumber \\              U_8(x,y) & =
   \frac{1-y}{y(x-1)}    && U_9(x,y)  =
   \frac{x(1-y)^2}{y(1-x)^2}       \nonumber   
\end{align}
We note ${\cal W}_{ {\cal S}{\cal K}}$ the planar web associated to
$U_1,..,U_9 $. Its singular locus is 
\begin{align} 
\Sigma_ { {\cal
    S}{\cal K}} = \{ (z,\zeta)\in  \rC ^2 \:
| \;  \;  z \; \zeta \, (1-z) \, (1-\zeta) (z-\zeta)&  (1+\zeta) (1+z)
    \times 
    \nonumber \\ 
     (1-z \zeta) -(  2-z-\zeta)& ( z \zeta -2\zeta +1)(2 z \zeta
-\zeta- z ) 
 =0  \; \}  \nonumber 
\end{align}
We choose again $\omega=(\frac{1}{3},\frac{1}{2} ) \in \rR ^2
\setminus \Sigma_ { {\cal
    S}{\cal K}} $ .
We want to find the local holomorphic solutions of $
({\cal SK}) $ at $\omega$.\\
As in the case of the dilogarithm, we get that $ U_i( \rC ^2 \setminus
\Sigma_ { {\cal S}{\cal K}} ) = \rC \setminus \{ 0,1 \} $  so, if 
$(F_1,..,F_9) \in  \underline{{\cal S}\scriptstyle{ \stackrel{\cal O}{{}\omega}}}
( {\cal SK}) $, then $ F_j \in \underline{ \widetilde{\cal
    O}\scriptscriptstyle{\omega} }    ( \rC \setminus \{0,1 \} ) $ for
$ j =1,...,9 $ .\\
In this case, the method of monodromy a priori can be applied to find
all the elements of $ \underline{{\cal S}{ \scriptstyle{\stackrel{\cal
        I}{{}_\omega}}}} ( {\cal SK}) $. Then by applying Abel's
method, we get next the missing solutions (noted $ {\bf F}_8, {\bf
  F}_{10}, {\bf F}_{15}, {\bf F}_{16}$ and $ {\bf F}_{17} $ below).\\

One can verify that the $28$
following  9-uplets of holomorphic germs are elements of 
$ \underline{{\cal S}{ \stackrel{\cal O}{{}_\omega}}}
( {\cal SK}) $ :

\begin{align}
&{\bf F}_1 =  \biggl( \, {\bf L}_{x_0}\, , \,- {\bf L}_{x_0}\, , \, -
 {\bf L}_{x_0}\,  ,0,0,0,0,0,0  \,  \biggr) \nonumber \\
&{\bf F}_2 =  \biggl( \, {\bf L}_{x_0+x_1}\,  ,\, 0\,  ,\,  - {\bf
 L}_{x_0+x_1}\,  ,\,  {\bf L}_{x_1}\, ,\, 0,  0,0,0,0  \,  \biggr) \nonumber  \\
& {\bf F}_3 =  \biggl( \, {\bf L}_{x_1} \, , \, {\bf L}_{x_1} \, ,\, 0
 \, ,
 -{\bf L}_{x_0}\, ,0,0,0,0,0 \,  \biggr) \nonumber \\
& {\bf F}_4 =  \biggl(\, 0,\, 0,\,  {\bf L}_{x_0}\, , {\bf L}_{x_0}\,,-
 {\bf L}_{x_0}\,,0,0,0,0 \,   \biggr) \nonumber   \\
& {\bf F}_5 =  \biggl( \,  {\bf L}_{x_1}\,  ,\,   0\, , \, - {\bf
  L}_{x_1} \, ,0\, , \, {\bf L}_{x_1} \, , 0,0,0,0 \,   \biggr)
 \nonumber \\
& {\bf F}_6 =  \biggl( \, {\bf L}_{x_0}\, , {\bf L}_{x_0}\,,0,0,0,- {\bf
  L}_{x_0}\, ,0,0,0,  \, \biggr)  \nonumber \\
& {\bf F}_7 =  \biggl( \,{\bf L}_{x_0}\, ,0,0, {\bf
  L}_{x_0}\,,0,0,-{\bf L}_{x_0}\, + i \pi ,0,0 \,   \biggr) \nonumber \\
& {\bf F}_8 =  \biggl( \, \iv \, , \, 0 ,\, 0,\, 0,\,  \iv  \,  ,\,
 0,\, \iv  -1\, , 0\, ,0 \, \biggr)  \nonumber \\
& {\bf F}_9 =  \biggl( \, {\bf L}_{x_1} \, ,0,0,0,0 , \, - {\bf
  L}_{x_1} \,  , \, {\bf L}_{x_1} \,  , 0\, , 0  \, \biggr)\nonumber \\
& {\bf F}_{10} =  \biggl( \, 0,\, \id  \, ,0,\, \id \,  ,\, 0,\, 0,\id -1, 0,0
\,  \biggr)  \nonumber      \\
& {\bf F}_{11} =  \biggl( \, 0 \, ,\, 0,\, 0,\, 0,\, 0, \, {\bf
 L}_{x_0}\,  , -{\bf L}_{x_0}\, ,  \,  {\bf L}_{x_0}\, ,0  \, \biggr)
  \nonumber \\
& {\bf F}_{12} =  \biggl( \,  0 \, , \, {\bf L}_{x_0}\, ,0\, ,0\, ,0\,
 ,0\, ,  {\bf
  L}_{x_1}\, , -{\bf L}_{x_1}\, ,0  \,  \biggr) \nonumber  \\
 & {\bf F}_{13} =  \biggl(\,  0,\,0,\, 0,\, 0,\, 0,\, 0,\, {\bf L}_{x_0}\, ,{\bf L}_{x_0}\,,-
 {\bf L}_{x_0}- 2 i {\pi}  \,  \biggr) \nonumber \\
& {\bf F}_{14} =  \biggl( \, 0,\, 0,\, 0,\, 0,\, {\bf L}_{x_1} \, ,\,
0, {\bf L}_{x_1} \, , \, 0,\, - {\bf L}_{x_1} \, \biggr) \nonumber \\
& {\bf F}_{15} =  \biggl( \, 0,\, \iv ,\ 0,\ 0,\, \id  , \, 0, \, 0,
\, \id -1 ,0 \,  \biggr)  \nonumber \\
& {\bf F}_{16} =  \biggl(  \, \id \, , \, 0 , \, 0, \iv \, ,\, 0,\,
0,\, 0,\, \iv  -1 , \, 0 \, \biggr)  \nonumber \\
& {\bf F}_{17} = \biggl( \, 0,\, 0,\, {\sf a}  ,\, 0, \, 0, \, -{\sf a} 
,\, 0,\, 0,\, -  {\sf a} \, \biggr) \nonumber \\
& {\bf F}_{18} =  \biggl( \, 2 \, {\bf L}_{x_0x_0 } \, ,\, 2 \, {\bf
  L}_{x_0x_0} \, ,\, -
 {\bf L}_{x_0 x_0  } \, , \, 0,\, 0,\, -{\bf L}_{x_0 x_0 } \, ,\, 0,\,
 0,\, 0  \,  \biggr) \nonumber  \\
& {\bf F}_{19} =  \biggl( \, 0,\, 0,\, 0,\, 0,\, 0, \, {\bf L}_{x_0x_0} \,
,-2 \, {\bf
  L}_{x_0x_0} \,,\, -2 \, {\bf L}_{x_0x_0} \, , \, {\bf L}_{x_0x_0}
\,+4i\pi\, {\bf L}_{x_0} \,-4\pi ^2  \, \biggr) \nonumber \\
& {\bf F}_{20} =  \biggl(\, 0,\, 0,  \, {\bf L}_{x_0x_0} \, ,-2 \, {\bf
  L}_{x_0x_0} \,, \, -2 \, {\bf L}_{x_0x_0} \,,\, 0,\, 0,\, 0,\, {\bf
  L}_{x_0x_0} \,  \biggr) \nonumber  \\
& {\bf F}_{21} =  \biggl( \, {\bf d} \, , \, -   {\bf d}\, , \, -
{\bf d}\, , \, -   {\bf d}\, ,   \, {\bf d}\, ,\, 
 0,\, 0,\,0,\, 0 \, \biggr) \nonumber  
\end{align}
\begin{align}
& {\bf F}_{22} =  \biggl( \,  {\bf d}\,  , \,   {\bf d}-
 \frac{i \pi }{2} \:  {\bf L}_{x_0} \, , \, 0, \, 0,\, 0, \, -
  {\bf d} \, , \,  {\bf d} \, ,  \, -  {\bf d}  \, , 0 \, \biggr)
 \nonumber \\
& {\bf F}_{23} =  \biggl( \, \pi ^2   ,0,0, {\bf d}   -\frac{i\pi}{2}{\bf
  L}_{x_0}, \, {\bf d} \, , \, 0,  \, {\bf d} \, , \, {\bf
  d}+\frac{i\pi}{2}{\bf L}_{x_0}+i\pi{\bf L}_{x_1}\,   , \, -  {\bf d}
\,   \biggr) \nonumber \\
& {\bf F}_{24} =  \biggl( \, {\bf L}_{x_0x_1}\,  ,\, {\bf
  L}_{x_0x_1}\,    , \, 0, \,  {\bf L}_{x_0x_0}\,,0,  \, -{\bf
  L}_{x_0x_1}\,  , \, {\bf L}_{x_0x_1}\,   ,    - \, {\bf
  L}_{x_0x_1}\,  
 - {\bf L}_{x_0x_0} + i \pi {\bf L}_{x_0} \, , \frac{\pi ^2 }{3}  \, 
 \biggr)  \nonumber \\
& {\bf F}_{25} =  \biggl(  \, 0,\,  {\bf L}_{x_0x_0} \, , \, 0 , \,
{\bf L}_{x_0x_1}\,,  \, {\bf L}_{x_0x_1}\,  , \, 0  , \,  {\bf
  L}_{x_0x_1} \,     , \,  {\bf L}_{x_0x_1} \,  , \,  -  {\bf
  L}_{x_0x_1} \,   \biggr)  \nonumber \\
& {\bf F}_{26} =   \biggl( \,  2   \,  {\bf L}_{x_0x_1}\,    , \, 0 
 ,\, - {\bf L}_{x_0x_1}\, , 0 , \, 2 \, {\bf L}_{x_0x_1}\,    , \, - {\bf
 L}_{x_0x_1}\,       , \, 2 \, {\bf L}_{x_0x_1}\,     ,  \, 0 , \, -
 {\bf L}_{x_0x_1}\, \biggr)  \nonumber \\
&{\bf F}_{27}= \biggl( \,  2 \, {\g} \, , \, 2{\g} \, ,
\, -{\g} \, , \, 2 \, {\g} \, , \, 2 \, {\g} \, , \, 
-{\g} \, , \, 2 \, \widehat{\g} \, 
, \,  2 \, \widehat{\g} \, , \, -{\g} \,  \biggr) \nonumber \\
  &{\bf F}_{28}= \biggl(  2\h(\bullet),2\h(\bullet)-\frac{2 \pi ^2}{3}{\bf
  L}_{x_0} ,
-\h(\bullet),2\h(\bullet),2\h(\bullet),
-\h(\bullet),2\widehat{\h}(\bullet)
,2\widehat{\h}( \bullet),- \h(\bullet) \biggr) \nonumber 
\end{align}
with
\begin{align} 
\id \; & := {\bf I} {\sf d}_{\mathbb C}  && \iv := {\bf 1}/{\bf I}{\sf d}_{\mathbb C} \nonumber \\
{\sf a} & : \bullet \rightarrow {\bf a}\mbox{rcth} \:  ( \sqrt{
  \bullet  }\: )
&& {\bf d}  := {\bf L}_{x_0x_1} - {\bf L}_{x_1 x_0} -\frac{\pi
  ^2}{6} \nonumber \\
{\g} &  :=2\, {\bf L}_{x_0^2 x_1} -{\bf L}_{x_0 x_1x_0}-{\bf L}_{
x_1x_0^2}- {{\frac{2}{3}\l{3} (1) }} &&
\widehat{\g}    := \g+{i\pi}{\bf L}_{x_0x_1}-4i\pi {\bf L }_{x_1x_0}
 -{\pi}^2 {\bf L}_{x_1}+2i{\pi}^3   \nonumber \\
 \h &  := {\bf L}_{x_0 ^2 x_1}- {\bf L}_{x_1x_0^2}
 && 
\widehat{\h} := \h -i \pi {\bf L}_{x_0x_1} 
   +2i \pi {\bf L}_{x_1x_0} +\frac{{\pi}^2}{2}
  {\bf L}_{x_1} -\frac{2i\pi ^3}{6} &&  \nonumber 
\end{align} 
Let $ \{ {\bf F}_l \, | \, l=-7,...,0 \, \} $ be a basis of the space
of the constant solutions of ${\cal SK}$~. Then it's just a tedious exercise of linear algebra to verify that the
$ {\bf F}_i$ 's ( for $ -7 \leq i \leq 28 $ ) are $36$ linearly independent
elements of $ \underline{{\cal S}{ \stackrel{\cal O}{{}_\omega}}}
( {\cal SK}) $ . Then it comes that
$$ 36 =     \mbox{ dim} _{ \scriptstyle{ \rC}}  \, \left< \, \{ \, {\bf F}_j
\} \, \right> \leq 
\mbox{ dim} _{ \scriptstyle{ \rC}} \, \underline{{\cal S}{ 
\stackrel{\cal O}{{}_\omega}}}   ( {\cal SK})   \leq {9(9-1)}/{2}=
36 $$
(the last inequality comes from proposition 4) .\\
So we have 
$$    \underline{{\cal S}{ 
\stackrel{\cal O}{{}_\omega}}}   ( {\cal SK})  = \left< \, \{ \, {\bf F}_j \,
| \, -7 \leq j \leq 28   \, \} \, \right> $$ 
This solves $ \eq{SK} $ at $ \omega$ in the holomorphic class.
We get the the local holomorphic solutions around $ \omega'\not \in
\Sigma_{\cal {\cal SK} } $ by analytic continuation of the $ {\bf F}_j$'s.        
\subsection{ An {\sf Afe} associated to a degenerate configuration of 5
  points }
Here we are considering the following {\sf Afe} 
\begin{align} 
 G_1(x) \;  +   \; & G_2(y) \; +  \; G_3(\frac{x}{y}) \;  +\;
          G_4(\frac{1-y}{1-x}) + \; 
          \;  G_5(
          \frac{x(1-y)}{y(1-x)}  \qquad \qquad  \qquad \qquad ( {\cal E}_{ {\sf
          c}}) \nonumber \\ 
+&  \; G_6(\frac{1+x}{1+y}) \; +   \; G_7(
\frac{x(1+y)}{y(1+x)}) \; + \; G_8(  \frac{(1-y)(1+x)}{(1-x)(1+y)}) \;
          =0 \nonumber
\end{align} 
We set $ V_6(x,y)=\frac{1+x}{1+y} $
,  $ V_7(x,y)=\frac{x(1+y)}{y(1+x)}$ ,  $
 V_8(x,y)=\frac{(1-y)(1+x)}{(1-x)(1+y)} $ and  $ V_i = U_i
$  for $i=1,..,5$~. \vspace{0.1cm} \\
We note $ {\cal W}_{ \sf c} $ the web associated
to the $V_i$'s: we will see in the next part that it is associated to
a configuration of 5 points. 
 A simple computation gives us its singular locus $ \Sigma_{ {\sf
    c}}$ . \\
We take $ \omega =(1/3,1/2) \in \rR \setminus \Sigma_{ {\sf c}} $.
 We want to determine the space $ \underline{{\cal S}{ 
\stackrel{\cal O}{{}_\omega}}}   ( {\cal E}_{ \sf c}) $.\\

By applying Abel's method, the author has constructed the following 21 elements
of $ \underline{{\cal S}{ 
\stackrel{\cal O}{{}_\omega}}}   ( {\cal E}_{ \sf c}) $\,:

\begin{align}
& {\bf G}_{1}=\biggl( \, 0\, ,\, 0\, , \, 
2 \, {\bf j} \, , \, {\bf j}\,  ,\,  0\, ,\, -{\bf j}\, ,\, 0,\, -1 \, 
\biggr) 
&& {\bf G}_{2}=\biggl(\,  0\, ,\, 0\, ,\, 2 \, {\bf j} \, 
,\, 0\, ,\, -{\bf j}\,  ,\, 0\, ,\, -{\bf j}\,  ,\, 0 \, 
\biggr) \nonumber \\
& {\bf G}_{3}=\biggl( 1,0,0,0,0,
{\bf j},-{\bf j}, \, -2 \, {\bf j}
\biggr)
&& {\bf G}_4=\biggl( {\bf L}_{x_0},- {\bf L}_{x_0},-
{\bf L}_{x_0},0,0,0,0,0 \biggr) \nonumber \\
& {\bf G}_5 =  \biggl( {\bf L}_{x_0+x_1} ,0,-{\bf
 L}_{x_0+x_1} ,  {\bf L}_{x_1} ,0,0,0,0  \biggr) 
  &&{\bf G}_6 =  \biggl( -{\bf L}_{x_1} , {\bf L}_{x_1},0,
 {\bf L}_{x_0} ,0,0,0,0  \biggr) \nonumber  \\
& {\bf G}_7 =  \biggl(0,0, {\bf L}_{x_0} , {\bf L}_{x_0},-
 {\bf L}_{x_0},0,0,0 \biggr) &&
{\bf G}_8 =  \biggl(  {\bf L}_{x_1} ,  0, -{\bf L}_{x_1} ,0, {\bf
  L}_{x_1},0,0,0  \biggr) \nonumber  \\
& {\bf G}_9 =  \biggl(  {\bf L}_{x_0 -x_{-1}}  ,  0,
 -{\bf L}_{x_0+x_{1}},0,0,{\bf
 L}_{x_0+x_1} ,0,0   \biggr) 
&& {\bf G}_{10} =  \biggl(  {\bf L}_{x_{-1}}  , -{\bf L}_{x_{-1}},
0,0,0,-{\bf L}_{x_{0}},0,0 \biggr) \nonumber\\
& {\bf G}_{11} =  \biggl( 0,0 ,   {\bf L}_{x_0} ,0,0, -{\bf L}_{x_0}, -{\bf L}_{x_0},0 \biggr) 
&& {\bf G}_{12} =  \biggl(  {\bf L}_{x_{-1}} , 0,{\bf L}_{x_1}, 0,0,0,-{\bf L}_{x_1},0\biggr) \nonumber \\
& {\bf G}_{13}=\biggl(0,0,0, {\bf L}_{x_0},0,{\bf L}_{x_0},0,
-{\bf L}_{x_0} \biggr) 
&& {\bf G}_{14} =  \biggl( -{\sf L}_2 ,  {\bf L}_{x_{-1}} ,0,
{\bf L}_{x_{1}}, 0,0,0,-{\bf L}_{x_{1}} \biggr) \nonumber 
\end{align}
\begin{align}
& {\bf G}_{15} =  \biggl( \, {\bf d}, \, -   {\bf d},\,   -
{\bf d} \, ,  \, -   {\bf d} \, ,   \, {\bf d} \, , \, 
 0 \,, \,0 \,, \,0 \,, \,0 \, \biggr)  \nonumber \\
&{\bf G}_{16} =  \biggl(  \, {\bf L}_{x_{-1} x_0 -x_0 x_{-1}}  \,  ,
\,  - {\bf L}_{x_{-1} x_0 -x_0 x_{-1}}  \,   , \,  - 2\, 
{\bf d} \,  ,  \,  0 \, , \,  0 \,   ,  \, 2\,   {\bf d} \,  ,  \,
 2\,   {\bf
  d}  \, ,  \,0 \biggr)  
\nonumber  \\
& {\bf G}_{17} =\biggl( \, 0 \,  , \,  0 \,  , \,   \tiny{8} \, {\bf c}- \tiny{10} \, {\bf j}, - {\bf c}-{\bf j} ,
 - {\bf c} ,  - {\bf c}-{\bf j}, - {\bf c}    +4 \, {\bf j} ,  8 \, {\bf c} 
\biggr) \nonumber  \\
& {\bf G}_{18} =  \biggl( \,  {\bf L}_{ x_{-1} x_1 -x_1 x_{-1} } +{\sf L}_2
{\bf L}_{ x_1 +x_{-1}}  \, , \,  -  {\bf L}_{ x_{-1} x_1 -x_1 x_{-1} }-
{\sf L}_2 {\bf L}_{ x_1 +x_{-1}}  \,    , \,  0 \,  ,  \,  - 2 {\bf
  d}\, , \,
0,  \,  - 2{\bf d} \,   , \,  0  \,  ,  \, 
 2\, {\bf d}  \,   \biggr)
\nonumber  \\
& {\bf G}_{19} =  \biggl(\, {\bf L}_{{x_0 x_1}+{x_0 x_{-1}}}\,  , \,  {\bf
  L}_{{x_1 x_0}+{x_{-1} x_0}} \, , \, 
 0 \, , \,  - {\bf L}_{ {x_1^2}  +  x_0 x_1} \, , \, {\bf L}_{ {x_1^2}  +
   x_0 x_1} \,  , \, {\bf L}_{ x_1^2+ x_1 x_0}\,  , \,   - {\bf L}_{ {x_1^2}  +
   x_0 x_1} \, ,  \,  -\frac{\pi ^2}{12} \,     \biggr) \nonumber  \\
& {\bf G}_{20} =\biggl( \, {\bf L}_{x_{-1}}\,  , \, - {\bf
  L}_{x_{-1}}\, , \, -2 \, {\bf b} \, 
, \, {\bf b}\, , \,   {\bf b}\,  , \,   {\bf b}\, , \, {\bf b}\,  ,\,
-2 \, {\bf b} \, \biggr) \nonumber \\
& {\bf G}_{21} =\biggl(    {\bf L}_{ x_{-1}x_1-x_0x_1},{\bf L}_{x_1 x_{-1}-x_1 x_0} 
 , {\bf L}_{  {x_1^2}  +  x_0 x_1 } -{\sf L}_2 {\bf L}_{x_0+x_1}  ,
 {\sf c}_{21} ,
\nonumber  \\
&  \qquad \qquad \qquad \qquad \qquad  -{\bf L}_{  {x_1^2}  +  x_0 x_1 } +{\sf L}_2 {\bf L}_{x_0+x_1}, 
 -{\bf L}_{  {x_1^2}  +  x_0 x_1 } +{\sf L}_2 {\bf L}_{x_0+x_1}, 0,
 {\bf L}_{x_1^2 +x_0 x_1} \biggr) \qquad \qquad  \nonumber  
\end{align}

with
\begin{align}
 {\bf j} & = \frac{1}{(1-\bullet)} &&
 {\bf b}   =\frac{{\bf
    L}_{x_0}}{1-\bullet}  &&& 
 {\bf c} = \frac{1}{(1-\bullet)^2}   \nonumber \\
{\sf L}_2 &  = \log(2) && 
 {\sf c}_{21}  =\frac{\pi ^2}{6}-\frac{1}{2} \log^2(2) &&&  \nonumber 
\end{align}
(the constant ${\sf c}_{21}$ had been determined numerically with the
precision of ${10}^3 $ digits).

Let $ \{ {\bf G}_l \, | \, l=-6,...,0 \, \} $ be a basis of the space
of the constant solutions of~$({\cal E}_{\sf c})$~.
Then it's just a tedious exercise of linear algebra to verify that the
$ {\bf G}_i$ 's ( for $ -6 \leq i \leq 21 $ ) are $28$ linearly independent
elements of $ \underline{{\cal S}{ \stackrel{\cal O}{{}_\omega}}}
( {\cal E}_{\sf c}   ) $ .\\
Then it comes that
$$ 28 =     \mbox{ dim}_{{}_{ \scriptstyle{ \rC}}}  \, \left< \, \{ \, {\bf G}_j
\} \, \right> \leq 
\mbox{ dim} _{ {}_{\scriptstyle{ \rC}}} \, \underline{{\cal S}{ 
\stackrel{\cal O}{{}_\omega}}}   ({\cal E}_{\sf c}   )   \leq {8(8-1)}/{2}=
28 $$
So we have 
$$    \underline{{\cal S}{ 
\stackrel{\cal O}{{}_\omega}}}   ({\cal E}_{\sf c} )  = \left< \, \{ \, {\bf G}_j \,
| \, -6 \leq j \leq 21  \, \} \, \right> $$ 
\section{Application to web theory and to the characterization of
  polylogarithmic functions of order $ \leq 3 $ 
  by their functional equation}
In the introduction we noticed that abelian functional equations arise
  in many areas of mathematics. We now give some 
applications of the material exposed in the preceding part to two of
  these areas: planar web
geometry and theory of polylogarithmic functional equations.

\subsection{Application to web theory}
 
\subsubsection{a brief introduction to planar web geometry}

We now briefly recall, in the analytic setting, the basic notions of
planar web geometry (the standard reference is  \cite{blabol}. See
\cite{akgo}, \cite{chern}, \cite{cherngriff1} or \cite{web} for more modern points of view).\\
A planar N-web ${\cal W}$ on a domain $ \Omega $ in a 2-dimensional
complex manifold $X$, is the data of an unordered set $ \{  {\cal F}_i
  \}$ of $N$ foliations of $\Omega$ such that their
  leaves are in general position. We are interested in the
  geometric local study of these webs to which we want to attach
  some invariants. A classical example of web is the ``algebraic N-web
  ${\cal W}_{ C}$'' associated to an algebraic reduced curve $ { C} \subset
  \projc^2$ of degree N: let $L_0$ be a generic line in $ \projc ^2$
  which transversally intersects
 the regular part of $C$ in N points\,: $ C.L_0=
P_1(L_0)+\cdots+P_N(L_0) $. There exists an open neighbourhood $
\Omega_{0} $ in the dual projective space $ \projc^{2 \star} $ 
and there are $N$ holomorphic maps $P_i \, : \Omega_{0} \rightarrow C$ 
such that all $ L \in \Omega_0$ transversally intersect  $C$ and  $
C.L=P_1(L)+\dots+P_N(L) $. Let $ {\cal F}_i $ be the foliation of $
\Omega_0$, the leaf of which at $L$ is the segment line $ \{ P_i=P_i(L)
\} $. Then $  {\cal W}_C=\{ {\cal F}_i \}_{i=1..N} $ is a web on $
\Omega_0$. Because the leaves of the foliations are segments of straight
lines, $  {\cal W}_C $ is a ``linear web''.

The problem of the linearization of (germs of) planar webs was a central
one. It has recently been solved in all its generality by M.Akivis,
V.Goldberg and V.Lychagin\,: for a
N-web, they give $N-2$ differential invariants which vanish if and
only if the web is linearizable (see \cite{akgoly}).

Let us consider our algebraic web $ {\cal W}_C $ again. Assume that
$C$ is smooth (to simplify) and let $\omega$ be a differential of the first kind on $C$. Abel's theorem implies that $ \sum {P_i} ^{\star} \omega=0$ on
$\Omega_0$: the abelian sums vanish. 

Let ${\cal W}=\{ {\cal F}_i \}_{i\leq N} $ be a (germ of) web at the
origin in $\rC ^2$. Then there exists $N$ germs of holomorphic maps
$U_i $   such that the leaves of $
{\cal F}_i $ are the level curves of $U_i$. We copy the notion of
abelian sum for this general web ${\cal W}$\,: an abelian relation for
${\cal W}$ (relatively to the $U_i$'s) will be an equation of the form $ \sum G_i(U_i)dU_i=0$ in
the space $ \underline{{\Omega}\scriptstyle{\stackrel{0}{1}}} $ of 
holomorphic germs of 1-form at the origin in $ \rC^2$. The space of
abelian relations of ${\cal W}$ relatively to the $U_i$'s will be noted
$ {\cal A}( \, {\cal W}) $. It has a natural structure of linear
space. By definition, the rank of $ {\cal W}$ will be $$ r_k( {\cal W}):= \mbox{dim}_{
  \scriptscriptstyle{\rC}}\, {\cal A}( \, {\cal W}) $$
The general
version of proposition 4 is that the rank is always finite and that 
$ r_k ( {\cal W} ) \leq (N-1)(N-2)/2$. The rank is a well defined  invariant for
webs (up to local diffeomorphisms).\\
Using the notations of this paper, we see that, modulo the constant
solutions, we have an isomorphism between 
$ \underline{{\cal
      S}{ \scriptstyle{\stackrel{\cal O}{\omega}}}} ( \, {\cal E}\{U_i\}
  )  $  and $ {\cal A}(  {\cal W} ) $ where $ {\cal E}\{U_i\} $ is the 
$ {\sf Afe}$ $ F_1(U_1)+\dots+F_n(U_N)=0 $. Then the results of part 2
may be seen, in the framework of web geometry, as tools to study the
abelian relations of webs the foliations of which are the level curves of
real rational functions. We will see that such webs are of particular interest.

For our algebraic web $ {\cal W}_C $, because 
$$  \mbox{dim}_{
  \scriptscriptstyle{\rC}}\, H^0(C,\Omega_C^1)=\frac{(N-1)(N-2)}{2}\geq r_k(
{\cal W}_C ) $$  
we obtain that the following map is an isomorphism:
\begin{align}
 H^0(C,\Omega_C^1) \, & \longrightarrow \, {\cal A}(\, {\cal W}_C)
 \nonumber \\
\omega^{\lambda}  \quad  & \longrightarrow \: \sum {P_i}^{\star} \omega
 ^{\lambda}   \nonumber
\end{align} 
Therefore an algebraic web is a linear web of maximal rank.\\
The Abel inverse theorem (due
to Lie for $N=4$, and generalized by Poincar{\'e}, Darboux, Griffiths
\cite{griff}  and Henkin \cite{henkin}) 
 tells us that a web of maximal rank is algebraic if and only
 if it is linearizable. 

It is well known that N-webs of maximal rank are linearizable (therefore
algebraic) when $N=3,4$ (the case $N=3$ is easy, and $N=4$ is due
to some work by Lie on translation surfaces. A naive idea would be
that all maximal rank web are linearizable and therefore algebraic. But this
is no longer true for $N$-webs with $N\geq 5$\,: in \cite{bol2},
G. Bol gave an example of a non-linearizable 5-web of maximal rank
which cannot be algebraic: this web (now known as `` Bol's
web'' and noted $ {\cal B}$) is the web the foliations of which are
the level curves of the $U_i$'s of equation $ ({\cal R})$ in part 3.3.
From the elements $ {\bf
  \Delta}_1,  {\bf \Delta}_2, {\bf \Delta}_3,  {\bf \Delta }_4, {\bf
  \Delta}_5 $ and  $   {\bf \Delta}_{6}$  of  $    \underline{{\cal
      S}{ \scriptstyle{\stackrel{\cal O}{\omega}}}}  ( \cal R)  $ 
 we can construct a basis
  of the space ${\cal A}(\cal B)$ of the abelian relations of $\cal B$
  at any generic point~$\omega_0$.

So we have $\mbox{dim}_{{\scriptscriptstyle{\mathbb C}}} \, {\cal A}({\cal B}) =6 $
  and  the web $\cal B$ is of maximal rank 6, although it is not 
linearizable.
From its discovery by Bol in the 30's onwards, this has been the single
  known counterexample 
  to the problem of linearization of planar webs of maximal rank. \\

Such webs, which looked very special, are called ``exceptional
webs'' (see section 6 in \cite{akgo} and part 3.2 and 3.3  of
\cite{hen1} for the problem of linearization of webs of maximal rank).
\subsubsection{exceptional planar webs and configuration of points}

According to Chern and Griffiths (see \cite{cherngriff2} page 83),
classifying the non linearizable maximal rank webs is the fundamental
problem in web geometry.  Since Bol's web is related to the
functional equation of Rogers dilogarithm, the ``Spence-Kummer web'' $
{\cal W}_{ \cal S K} $ associated to the Spence-Kummer equation of the
trilogarithm ``{\it seems to be a good candidate as an
  exceptional 9-web} '', as noticed by A. H{\'e}naut in part 3.3 of
\cite{hen1}. We now prove that this web actually is exceptional. The explicit resolution of the equations $
\eq{SK} $ and $ ( {\cal E}_{\sf c})$ done respectively in parts 3.4 and
3.5 allows us to find other new examples of such
``exceptional webs''.\\

We first study $ {\cal W}_{\cal SK} $ and its subwebs.\\
 For any subset $J\subset \{ 1,..,9\}$ we note ${\cal W}_J$
  the $|J \, |$-subweb of $ {\cal W}_{\cal SK} $ given by the level curves of the
  function $U_j$, with $j\in J$. 
If $ j_1,..,j_p$ are p distinct integers in $ \{1,..,9\}$ , then we
  note $ \widehat{j_1..j_p}:=\{1,..,9\} \setminus \{ j_1,..,j_p\}$.
\begin{theorem}{\bf :}  \newline
\begin{tabular}{ll} 
 \qquad $ \bullet$ $ {\cal W}_{\cal SK} $  is an exceptional 9-web &  $ \bullet$  ${\cal W}_{\widehat{ 69}} $
 is an exceptional 7-web \\
 \qquad$ \bullet$  ${\cal W}_{\widehat{ 679}} $  is an exceptional 6-web &
   $ \bullet$  ${\cal W}_{\widehat{ 248}} $  is an exceptional 6-web \\
 Those two exceptional 6-webs are not equivalent. &   \\
 \qquad  $  \bullet$ ${\cal W}_{\widehat{ 369}} $  is an
 hexagonal 6-web. &   
\end{tabular}\\
\end{theorem}
{\bf proof:} 
For each of these webs, we have to prove two distinct things: the first
is that the rank is maximal, the second is that the web is
non-linearizable.\\
From the basis $ \{ \, {\bf F}_j \} $ of $ {\underline{{\cal
      S}{ \stackrel{\cal O}{{}_\omega}}}   ( {\cal W}_{\cal SK}) } $
  described in 3.4, we can easily construct $28$ linearly
  independent abelian equations for $ {\cal W}_{\cal SK} $, which is of maximal rank. For $\omega$ generic, we have
 natural linear inclusions $   \underline{{\cal
      S}{ \stackrel{\cal O}{{}_\omega}}}   ( {\cal W}_J)
  \hookrightarrow   \underline{{\cal
      S}{ \stackrel{\cal O}{{}_\omega}}}   ( {\cal W}_{\cal SK}) $.
 From this we can easily deduce an explicit basis of the spaces  
$ \underline{{\cal S}{ \stackrel{\cal O}{{}_\omega}}}   ( {\cal
  W}_J)$, and so of the space $ {\cal A} \, ( {\cal W}_J)$ for any subset
$J \subset \{1,..,N\} $. So we can calculate the rank of any sub-web
of $ {\cal W}_{\cal SK} $: all the webs in proposition 4 have maximal
rank.\\
Let us note $ {\cal W} $ an exceptional web of proposition 6 distinct from $ {\cal
  W}_{ \widehat{248}}$. We remark that ${\cal W}$ contains ${\cal B}$
  as a 5-subweb. Because ${\cal B} $ cannot be linearized, the same is
  true for $ {\cal W} $ which thus is exceptional. We can't use this argument to prove that $ {\cal
  W}_{ \widehat{248}}$ is not linearizable: it's easy (but tedious) to
  see that all the 5-subwebs of ${\cal W}_{\widehat{ 248}} $ have rank
  5 and so are not equivalent to Bol's web (this already shows 
  that ${\cal W}_{\widehat{ 248}} $ and ${\cal W}_{\widehat{ 679}} $
  are not equivalent).\\
One can verify that its
  associated polynomial $ P_{ {\cal
  W}_{ \widehat{248}}   } $ (see \cite{henaut} for a definition) is
  of degree $4>3$. 
Theorem 2 in \cite{henaut} says that if $ \widetilde{\cal W} $ is a
web of maximal rank, then it is linearizable if and only if $ 
P_{ \widetilde{\cal W}   } $ is of degree smaller than 3.
Because $ {\cal  W}_{ \widehat{248}}$  is of maximal rank, it implies
that it is exceptional. $ \blacksquare$ \vspace{0.5cm}  \\
{\bf remarks :} 
{\bf 1.}  The sub-webs ${\cal W}_{\widehat{ 36}} $ and 
${\cal W}_{\widehat{ 39}} $ are exceptional too but equivalent to
  ${\cal W}_{\widehat{ 69}} $~. \\
{\bf 2.} The sub-webs    ${\cal W}_{\widehat{ 689}} $, 
${\cal W}_{\widehat{ 349}} $,
${\cal W}_{\widehat{ 236}} $,
${\cal W}_{\widehat{ 359}} $, and
${\cal W}_{\widehat{ 136}} \, $
 are exceptional too
  but equivalent to ${\cal W}_{\widehat{ 679}} \, $.\\
 {\bf 3.} The sub-webs   ${\cal W}_{\widehat{ 147 }} $ ,
${\cal W}_{\widehat{ 257}} $ , and 
${\cal W}_{\widehat{ 158 }} \, $ 
are exceptional too
  but equivalent to ${\cal W}_{\widehat{248}} \, $.\\
{\bf 4.} The exceptional d-subwebs of $ {\cal W}_{\cal SK} \, $(with
$ d \geq 6$ ) are those which are described in proposition 9 and in the
above remarks {\bf 1},{\bf 2} and  {\bf 3} .\\
{\bf 5.} We have a beautiful functional equation associated to
  ${\cal W}_{\widehat{248}} $ for $ \l{2} $. It is given by the element ${\bf
  F}_{26}$ of part 3.4 : 
\begin{center} 
$ 2\l{2} (x) -\l{2}(\frac{x}{y})+ 2\l{2} (\frac{x(1-y)}{y(1-x)})
-\l{2}({x}{y}) + 2\l{2}
(-\frac{x(1-y)}{1-x})-\l{2}(\frac{x(1-y)^2}{y(1-x)^2 }) =0$ 
\end{center}
It is equivalent to Newman's functional equation of the bilogarithm (see formula
  (1.43) in \cite{lewin}, page 13).\\

By Bol's theorem (see \cite{blabol} page 108), the fact that ${\cal
  W}_{\widehat{ 369}}$ is hexagonal implies that it is linearizable into
  a web formed by 6 pencils of lines, therefore, by duality,  it is associated to a
  configuration of $6$ points on $ \projc ^2$. A linearization for
  this web is given by the quadratic Cremona transform $ {\bf {\sf C}} : \, (x,y)
  \rightarrow (1/(x-1),1/(y-1)) $ . 
It is natural to ask what is the action of $ {\bf {\sf C}} $ on the whole web 
$ {\cal W}_{\cal SK}$ .
We introduce some definitions. For $d>0$, let $\delta_d$ be the dimension of
the space of algebraic curves of degree $d$ in $ \projc ^2$. 
\begin{definition}
{\em Let be $ d\geq 1$. If
$ K$  is a set of $\delta_d-1$ points in general position in the complex
projective plane,
then the family  of the curves of degree
$d$ through these $\delta_d-1$ points is noted ${\cal F}_{K} $
. It is a singular foliation of $\projc^2$.}
\end{definition}
For $n\geq 3$, we note $ \Delta_n=\bigcup_{i<j} \{ (p_1,..,p_n))
\in {\bigl( \projc^2\bigr)}^n \, | p_i=p_j \} $.
We define the space of configuration of $n$ points in $ \projc^2$ as the set 
$ \underline{C}^n_2= {( \projc^2\bigr)}^n  - \Delta_n$. If three distinct points ${\sf p_i},{\sf
    p_j},{\sf p}_k$ of a configuration $({\sf p_1},..,{\sf p}_n)$ lie
on a same line, the configuration is said ``degenerate''. 
\begin{definition}
{\em Let be $N\geq 3$. The web $ {\cal W}_{\bf p} $ associated to a
  configuration $ {\bf p}=({\sf p_1},..,{\sf p}_N)$ of $N$ points in
  $\projc^2$ is the singular web defined on the whole plane, the
  foliations of which  are the $ {\cal F}_{J} $'s where $J$ runs on
  the set of subsets of $\delta_j-1$ points in $ \{ \, {\sf p
  }_1,...,{\sf p}_N \}
  $, in general position, with $  1\ \leq j \leq N $~.}
\end{definition}
It is well known that Bol's web ${\cal B}$ is associated to a
configuration of $4$ points in generic position in the projective plane.
More precisely, the web given by the level curves of the functions
$U_1,..,U_5$ is the web associated (in the sense of definition 4) to
the configuration $ { \sf \mathbf b}$ described by figure 1 below.\\
\begin{tabular}{ccc}
 \qquad  \qquad   \qquad  \qquad\begin{picture}(0,0)%
\includegraphics{nyobol.pstex}%
\end{picture}%
\setlength{\unitlength}{3947sp}%
\begingroup\makeatletter\ifx\SetFigFont\undefined%
\gdef\SetFigFont#1#2#3#4#5{%
  \reset@font\fontsize{#1}{#2pt}%
  \fontfamily{#3}\fontseries{#4}\fontshape{#5}%
  \selectfont}%
\fi\endgroup%
\begin{picture}(1573,1719)(949,-913)
\put(1730, 78){\makebox(0,0)[lb]{\smash{\SetFigFont{12}{14.4}{\rmdefault}{\mddefault}{\updefault}{\color[rgb]{0,0,0}$\projc ^1_{\infty}$ }%
}}}
\put(2409,-489){\makebox(0,0)[lb]{\smash{\SetFigFont{12}{14.4}{\rmdefault}{\mddefault}{\updefault}{\color[rgb]{0,0,0}$ {\sf b_1} $ }%
}}}
\put(1316,611){\makebox(0,0)[lb]{\smash{\SetFigFont{12}{14.4}{\rmdefault}{\mddefault}{\updefault}{\color[rgb]{0,0,0}$ {\sf b_2} $ }%
}}}
\put(1606,-368){\makebox(0,0)[lb]{\smash{\SetFigFont{12}{14.4}{\rmdefault}{\mddefault}{\updefault}{\color[rgb]{0,0,0}$ {\sf b_3} $ }%
}}}
\put(975,-529){\makebox(0,0)[lb]{\smash{\SetFigFont{12}{14.4}{\rmdefault}{\mddefault}{\updefault}{\color[rgb]{0,0,0}$ {\sf b_4} $ }%
}}}
\end{picture}

& \qquad \qquad \qquad \qquad \qquad & 
\begin{picture}(0,0)%
\includegraphics{nbonskum.pstex}%
\end{picture}%
\setlength{\unitlength}{3947sp}%
\begingroup\makeatletter\ifx\SetFigFont\undefined%
\gdef\SetFigFont#1#2#3#4#5{%
  \reset@font\fontsize{#1}{#2pt}%
  \fontfamily{#3}\fontseries{#4}\fontshape{#5}%
  \selectfont}%
\fi\endgroup%
\begin{picture}(1847,1966)(637,-1199)
\put(1730, 78){\makebox(0,0)[lb]{\smash{\SetFigFont{12}{14.4}{\rmdefault}{\mddefault}{\updefault}{\color[rgb]{0,0,0}$\projc ^1_{\infty}$ }%
}}}
\put(1316,611){\makebox(0,0)[lb]{\smash{\SetFigFont{12}{14.4}{\rmdefault}{\mddefault}{\updefault}{\color[rgb]{0,0,0}$ {\sf q_2} $ }%
}}}
\put(2409,-489){\makebox(0,0)[lb]{\smash{\SetFigFont{12}{14.4}{\rmdefault}{\mddefault}{\updefault}{\color[rgb]{0,0,0}$ {\sf q_1} $ }%
}}}
\put(1356,-1001){\makebox(0,0)[lb]{\smash{\SetFigFont{12}{14.4}{\rmdefault}{\mddefault}{\updefault}{\color[rgb]{0,0,0}$ {\sf q_6} $ }%
}}}
\put(637,-1141){\makebox(0,0)[lb]{\smash{\SetFigFont{12}{14.4}{\rmdefault}{\mddefault}{\updefault}{\color[rgb]{0,0,0}$ {\sf q_3} $ }%
}}}
\put(683,-526){\makebox(0,0)[lb]{\smash{\SetFigFont{12}{14.4}{\rmdefault}{\mddefault}{\updefault}{\color[rgb]{0,0,0}$ {\sf q_5} $ }%
}}}
\put(1349,-499){\makebox(0,0)[lb]{\smash{\SetFigFont{12}{14.4}{\rmdefault}{\mddefault}{\updefault}{\color[rgb]{0,0,0}$ {\sf q_4} $ }%
}}}
\end{picture}
 \\
Figure 1: & \qquad & Figure 2: \\
configuration ${\bf {\sf b}}$ with  & \quad & configuration ${\bf {\sf q}} $ with         \\
  $ \scriptstyle{ {\bf {\sf b_1=[1:0:0] , \, b_2=[0 :1:0]       }}}$ & &
$ \scriptstyle{ {\bf {\sf q_1=b_1 , \, q_2=b_2 , \, q_3=[-1:-1:1]    }}} $  \\
$  \scriptstyle{ {\bf {\sf  b_3=[1 : 1:1] }} ,  \,  {\bf {\sf
      b_4=[0:0:1] }}} $  &  &  $ \scriptstyle{ {\bf {\sf q_4=b_4 
        , \, q_5=[-1 :0:1] , \,    q_6=[0:-1:1 ]   }}}$
\end{tabular}  \\

For the Spence-Kummer web ${\cal W}_{\cal SK} \, $, we have the following
\begin{prop}
  The web $ {\cal W}_{{\bf q}} $ associated to the degenerate
  configuration $ {{\bf q}}$ of $6$ points in $ \projc ^2$ given by
  figure $2$ above is the image of $ \, {\cal W}_{\cal SK} $ by  ${\bf {\sf C}}$
\end{prop}
The web $ {\cal W}_{\sf c} $ in 3.5, which is of maximal rank 21, is
also associated to a configuration noted ${\sf c}$ and defined by figure 3
below.\\

Because configuration $ {\sf b}$ is a subconfiguration of
configuration $ {\sf c}$, Bol's web is a sub-web of $ {\cal
  W}_{\sf c} $. Then this web is non-linearizable and since it is of
maximal rank (see part 3.5), it comes the 
\begin{prop}
The web $ {\cal W}_{\sf c} $ associated to $ {\sf
  c}$ is an exceptional planar 8-web.
\end{prop}
\begin{tabular}{ccc}
 \qquad \qquad  \qquad \qquad   \begin{picture}(0,0)%
\includegraphics{nconf8.pstex}%
\end{picture}%
\setlength{\unitlength}{3947sp}%
\begingroup\makeatletter\ifx\SetFigFont\undefined%
\gdef\SetFigFont#1#2#3#4#5{%
  \reset@font\fontsize{#1}{#2pt}%
  \fontfamily{#3}\fontseries{#4}\fontshape{#5}%
  \selectfont}%
\fi\endgroup%
\begin{picture}(1825,1978)(659,-1172)
\put(1730, 78){\makebox(0,0)[lb]{\smash{\SetFigFont{12}{14.4}{\rmdefault}{\mddefault}{\updefault}{\color[rgb]{0,0,0}$\projc ^1_{\infty}$ }%
}}}
\put(1316,611){\makebox(0,0)[lb]{\smash{\SetFigFont{12}{14.4}{\rmdefault}{\mddefault}{\updefault}{\color[rgb]{0,0,0}$ {\sf c_2} $ }%
}}}
\put(2409,-489){\makebox(0,0)[lb]{\smash{\SetFigFont{12}{14.4}{\rmdefault}{\mddefault}{\updefault}{\color[rgb]{0,0,0}$ {\sf c_1} $ }%
}}}
\put(975,-529){\makebox(0,0)[lb]{\smash{\SetFigFont{12}{14.4}{\rmdefault}{\mddefault}{\updefault}{\color[rgb]{0,0,0}$ {\sf c_4} $ }%
}}}
\put(659,-931){\makebox(0,0)[lb]{\smash{\SetFigFont{12}{14.4}{\rmdefault}{\mddefault}{\updefault}{\color[rgb]{0,0,0}$ {\sf c_5} $ }%
}}}
\put(1606,-368){\makebox(0,0)[lb]{\smash{\SetFigFont{12}{14.4}{\rmdefault}{\mddefault}{\updefault}{\color[rgb]{0,0,0}$ {\sf c_3} $ }%
}}}
\end{picture}
& \qquad \qquad \qquad & 
\begin{tabular}{c} 
Figure 3: \\
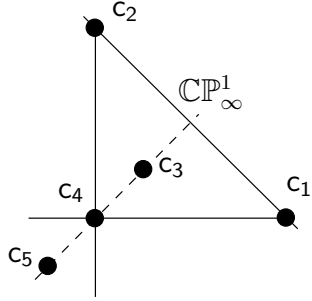
   configuration $ {\sf {\bf c}} $ which is \\
 associated to the web 
$ {\cal W}_{\sf c}$ of 3.5\\ with $ {\sf {\bf c_1}}={\sf {\bf b_1}}
,    {\sf {\bf c_2}}={\sf {\bf b_2}}    $           \\ 
$ {\sf {\bf c_3}}={\sf {\bf b_3}} , {\sf {\bf c_4}}={\sf {\bf b_4}} ,
 $ and  $ {\sf {\bf c_5}}={\sf {\bf q_3}} $
\\ \\ \\  
\end{tabular} 
\end{tabular}\\

The exceptional 6-subweb $ {\cal W}_{\widehat{679}}$   of $ {\cal
  W}_{\cal SK} $ is associated too with a
sub-configuration of~${\bf q}$: 
\begin{prop} The image of the exceptional web $ {\cal
    W}_{\widehat{679}}$ by $ {\sf C} $ is the web associated to the
  subconfiguration $ {\sf (q_1,q_2,q_3,q_4,q_5)} $ of $ {\bf {\sf
  q}}$.
\end{prop}
The other exceptional subwebs of $ {\cal W}_{\cal SK  }$ must also be 
associated to configurations of points but in a more complicated way
than in definition 4.\\

The fact that the only known exceptional planar webs described above are
related to configurations of points may be an important fact which
should be studied.\\

In \cite{damiano}, D. Damiano considers some webs of curves in $ \rR^N
 ,  \; (N\geq 2) $ similarly associated to configurations of
 points. He shows that those webs are exceptional curvilinear~webs.\vspace{0.2cm}\\
All this results allow to think that it could exist a real link between
configurations of points and exceptional webs.
In this spirit we have the following general results\,:
\begin{prop} The web associated to any configuration of 4 points in $
  \projc ^2$ is of maximal rank. It is non linearizable only if the
  configuration is generic\,: then it is (projectively) equivalent to Bol's web
  ${\cal B}$.
\end{prop} 
and for configurations of 5 points in the plane:
\begin{prop} The web associated to any degenerate configuration of 5
  points in $ \projc ^2$ is of maximal rank. 
\end{prop}
{\bf sketch of the proof\,:} 
We consider the stratification of $\underline{C}_5^2$ described by
figure 4: 
\begin{itemize}
\item ${\bf S}_0$ is the open subset of generic configurations.
\item ${\bf S}_1$ is the analytic strata of degenerate configurations such
  that three and only three points are lying on a same line .
\item ${\bf S}_2$ is the analytic strata of degenerate configurations such
  that exactly four points lie on a same line
\item ${\bf S}_3  $ is the analytic strata of degenerate
  configurations $({\sf p}_1,..,{\sf p}_5)$ outside $ {\bf S}_2$ such that there exists a
  unique $ {\sf p}_j $ such that for all $ i\neq j $ there exists $
  k$ distinct from $ i$ and $j$ such that the three points  $ {\sf p}_i,
  {\sf p}_j$ and $ {\sf p}_k $ are aligned.
\item ${\bf S}_4$ is the analytic strata of degenerate configurations such
  that the five points are aligned.

\end{itemize}
Each strata ${\bf S}_i$ is a
smooth connected analytic subvariety of~$\underline{C}_5^2$.\\
We note $N_0=10, \, N_1=8, \, N_2=5, \, N_3=6, \mbox{ and } N_4=5$. For each $i
\in \{0,..,4\}$, the web $ {\cal W}_{\sf p}$ associated to a
configuration $ {\sf p} \in {\bf S}_i$ is a $N_i$-web.\\

\begin{tabular}{c}
 Figure 4: \\
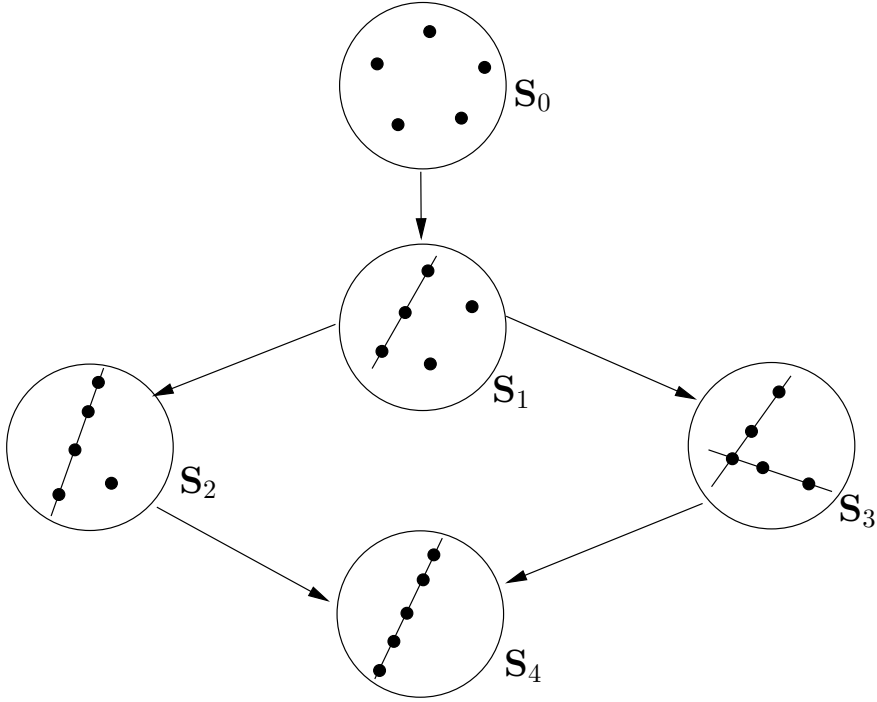
 stratification of $\underline{C}_5^2$ by degenerate configurations \\
  An arrow $ {\sf A} \rightarrow {\sf B} $ between two stratas  
  means that $  {\sf B} \subset {\partial} {\sf A} $ in
 $\underline{C}_5^2$ \\  \\
\qquad  \qquad \begin{picture}(0,0)%
\includegraphics{test.pstex}%
\end{picture}%
\setlength{\unitlength}{3947sp}%
\begingroup\makeatletter\ifx\SetFigFont\undefined%
\gdef\SetFigFont#1#2#3#4#5{%
  \reset@font\fontsize{#1}{#2pt}%
  \fontfamily{#3}\fontseries{#4}\fontshape{#5}%
  \selectfont}%
\fi\endgroup%
\begin{picture}(5343,4383)(456,-3514)
\put(1486,-2206){\makebox(0,0)[lb]{\smash{\SetFigFont{14}{16.8}{\rmdefault}{\mddefault}{\updefault}{\color[rgb]{0,0,0}$ \; {\bf S}_2$}%
}}}
\put(3451,-1636){\makebox(0,0)[lb]{\smash{\SetFigFont{14}{16.8}{\rmdefault}{\mddefault}{\updefault}{\color[rgb]{0,0,0}$ \; {\bf S}_1$}%
}}}
\put(5626,-2386){\makebox(0,0)[lb]{\smash{\SetFigFont{14}{16.8}{\rmdefault}{\mddefault}{\updefault}{\color[rgb]{0,0,0}$ \; {\bf S}_3$}%
}}}
\put(3526,-3361){\makebox(0,0)[lb]{\smash{\SetFigFont{14}{16.8}{\rmdefault}{\mddefault}{\updefault}{\color[rgb]{0,0,0}$ \; {\bf S}_4$}%
}}}
\put(3652,226){\makebox(0,0)[lb]{\smash{\SetFigFont{14}{16.8}{\rmdefault}{\mddefault}{\updefault}{\color[rgb]{0,0,0}$ {\bf S}_0$}%
}}}
\end{picture}
 
\end{tabular}\\

The natural action of  $ PGL_3( \rC)$ on $ \projc^2$ induces a group
action  $ {\bf {\sf q}}=({\sf q}_1,.., {\sf q}_5)
\rightarrow  {\bf {\sf q}}^g=(g({\sf q}_1),.., g({\sf q}_5))$
 on $\underline{C}_5^2$. Two webs $ {\cal W}_{ {\bf {\sf q}}}$ and $ {\cal
  W}_{ {\bf {\sf p}}}$ are projectively equivalent if and only if $ { {\bf {\sf
      q}}}$  and $ { {\bf {\sf p} }} $ belong to the same orbit. Then for
any orbit $ {\sf O} \subset \underline{C}_5^2 $, we consider a
particular configuration $ {\bf {\sf
    p}}_{\sf O} \in {\sf O} $. We prove that the rank of $ {\cal W}_{
  {\bf {\sf p}}_{\sf O}  } $ is maximal by constructing a basis of the
space $ {\underline{{\cal
      S}{ \stackrel{\cal O}{{}_\omega}}}   ( {\cal W}_{{\bf {\sf
      p}}_{\sf O} }  )}$ at a generic $\omega$. Moreover the rank is a
local invariant of the webs. All this implies proposition 10. We skip
here the explicit determinations of the spaces of abelian relations of
the webs $ {\cal W}_{
  {\bf {\sf p}}_{\sf O}  } $.   $ \blacksquare$\\

{\bf remarks:} \begin{enumerate}
\item For $a \in \rC \setminus\{0,1\}$, we note $ {\sf c}^a_5=[a:a:-1] $
  and $ {\sf c}_a=({\sf c}_1,..,{\sf c}_4,{\sf c}^a_5) \in {\bf S}_1$.\\
  The web $ 
{\cal W}_{{\sf c}_a} $ is exceptional: it is non-linearizable and, as
in part 3.5, we can construct a basis of dimension 21 of its space of
abelian relations.  This gives us
 a family of exceptional webs non-projectively equivalent. It would be
 interesting to know if they are locally equivalent or not.

\item For a generic configuration of 5 points, this  proposition is not
proved for the moment.
\end{enumerate}
All the results above allow us to state the 
\begin{conjecture}
The web $ {\cal W}_{{\bf q}} $ associated to any configuration ${\bf
  q}$ of points in the projective plane is of maximal rank.
\end{conjecture}
By an argument used above, we see that $ {\cal W}_{{\bf q}} $ is non linearizable as soon as ${\bf q}$
  contains a sub-configuration of $4$ points in general position, 
  so conjecture 1 may give us a list of exceptional webs.

It is under the inspiration of this conjecture than the author has
studied the {\sf Afe} $ ({\cal E}_{\sf c})$.\\

The results of part 3 show that most abelian functional
equations of webs associated to configurations of points studied in
part 4.1 are constructed
from iterated integral functions. If conjecture 1 is true, there could 
exist numerous {\sf Afe} linked to the
exceptional webs associated to configurations. We can expect than
some of those {\sf Afe} may be constructed from iterated
integrals too. This could be a way to find new functional equations for
higher order polylogarithms, which would be useful for the
K-theoretical study of algebraic number fields (see for instance \cite{zagier} or \cite{gangl}).

\subsection{application to the problem of characterizing polylogarithmic
  functions by their functional equation}
Our objective here is to study the function which satisfies the equation
$ (L_2)$ or $ ( SK)$. \\
This kind of problem has been studied for a long time for the Cauchy
equation $({\cal C})$: we know that any non-constant measurable
local solution of $({\cal C})$ is constructed from the logarithm.\\
The explicit resolution of equations $ \eq{R} $ and $ ( {\cal SK})$
done in part 3 allows us to get the same kind of results for the
 dilogarithm and the trilogarithm: these functions are respectively
 ``characterized'' (in the measurable class) by their functional equation $ (L_2) $ and $ (SK)$.
 \subsubsection{Characterization of the dilogarithm by Rogers
   equation $ ({\cal R})$}
We first have this result which comes easily from a result established
by Rogers in the early 20th century
(see \cite{roger} section 4) 
\begin{prop}  
If $F$ is a function  of class $C^3$ on $]\, 0,1[$ such that 
$$ ( \ast)  \qquad F(x)-F(y)-F(\frac{x}{y})-F( \frac{1-y}{1-x})+F(
\frac{x(1-y)}{y(1-x)})= 0 \qquad $$
for $ 0<x<y<1 $,  then we have $ f=  \alpha \, {\bf d} $ where  $ \alpha \in \rR$.
\end{prop}
The proof is essentially an application of Abel's method to this
 case. Proposition 1 allows us to see that proposition 11
 is still valid under the weaker assumption of
 measurability on $F$. \\
The dilogarithm has a single valued version: the Bloch-Wigner function 
$$ {\cal L}_2(z)= \Im m \biggl( \l{2} (z)  + \log(1-z) \log|z|  \biggr)
$$
which is real analytic on $ \projc^1 \setminus \{0,1,\infty\}$ and 
extends to
the whole projective line by continuity. 
For this function, the functional equation $ (L_2)$ becomes 
$$ (\Diamond) \qquad \qquad  \sum_{i=1}^4 (-1)^i {\cal
  L}_2(c_r(z_0,..,\hat{z_i},..,z_4))=0 \qquad , \quad 
z_i \in  \projc^1 \qquad \qquad    $$ 
where $ c_r(z_1,...,z_4)$ denotes the cross ratio of $4$
points. The equation $(\Diamond)$ takes the form $ (\ast)$ when we take
$(z_1,..,z_4)= (\infty,0,1,y,x)$.\\
In \cite{bloch}, Bloch proves the following characterization of $ {\cal
  L}_2$ in the measurable class by the equation $(\Diamond) $. 
\begin{prop}
Let $f : \, \projc ^1    \rightarrow \rR$ be measurable and satisfying~$(\Diamond) $. Then $f$ is proportional to $ {\cal L}_2$.
\end{prop}
Using proposition 1 and Bol's discovery of the space $ {\cal A}(
{\cal B})$ (see 3.3), we can state 
\begin{prop}
If $F,G$  are measurable functions on $]0,1[$ satisfying 
$$  F(x)-F(y)-F(\frac{x}{y})-F( \frac{1-y}{1-x})+G(
\frac{x(1-y)}{y(1-x)})= 0 $$ for $ 0<x<y<1$ 
,  then we have $ F=G=  \alpha \, {\bf d} $ with $ \alpha \in \rR$.
\end{prop}
(We can prove this result by a direct application of Abel's
method, because by proposition 1, the functions $F$ and $G$ of
proposition 14 are analytic). In the class of measurable functions, this result gives us a semi-local characterization of Roger's dilogarithm by its functional
equation $(R)$, for two unknown functions. In a certain sense, it's stronger than the results by
Rogers and Bloch. The explicit knowledge of $ \underline{{\cal S}{ 
\stackrel{\cal O}{{}_\omega}}}   ({\cal R} ) $ allows us to state
numerous variants of proposition~13.\\

Those results can be formulated in an inhomogeneous form to obtain
some characterization of $ \l{2} $ by functional equations inspired from $(L_2)$.

\subsubsection{characterization of the trilogarithm by Spence-Kummer
  equation $ ({\cal SK})$.}
The fact that the logarithm and dilogarithm are characterized by the   
{\sf Afe} with rational inner functions which they verify  naturally
leads us to ask if the same is true for any  trilogarithmic function.\\

In his paper \cite{gon}, A. Goncharov obtains some results of
this kind: 
he considers the real single-valued cousin of $ \l{3} $  introduced by
Ramakhrishnan and Zagier :
$$  {\cal L}_3(z):=
 \Re e\left( \l{3} (z)-\l{2} (z) \log|z| +\frac{1}{3}\l{1} (z) \log|z|^2 
 \right) $$ 
defined on the whole $\mathbb C \mathbb P ^1$ and extended 
to $ \mathbb R [ \mathbb C \mathbb P ^1] $ by 
 linearity . \\
When it is well defined, he considers the following element of 
$ \mathbb Q [ \mathbb C \mathbb P ^1] $ :
\begin{align}
 R_3( \alpha_1, \alpha_2,\alpha_3):= \sum_{i=1}^3
& \biggl(  \{ \alpha_{i+2} \alpha_i -\alpha{i}+1 \} +   
             \{ \frac{\alpha_{i+2} \alpha_i -\alpha_i+1}{\alpha_{i+2}
               \alpha_i} \} +\{ \alpha_{i+2} \} \nonumber \\
&  + \{ \frac{\alpha_{i+2} \alpha_{i+1} -\alpha_{i+2}+1}{
   (\alpha_{i+2} \alpha_i -\alpha_{i}+1)    \alpha_{i+1} } \}  -     \{ \frac{\alpha_{i+2} \alpha_i
   -\alpha_i+1}{\alpha_{i+2}} \}     - \{ 1 \}
  \nonumber \\
&    - \{ \frac{\alpha_{i+2} \alpha_{i+1} -\alpha_{i+1}+1}{ (\alpha_{i+2} \alpha_i
   -\alpha_i+1) \alpha_{i+1} \alpha_{i+2}  } \}  + \{- \frac{\alpha_{i+2} \alpha_{i+1} -\alpha_{i+1}+1)\alpha_i}{ \alpha_{i+2} \alpha_i
   -\alpha_i+1 } \}     \biggr)    \nonumber \\             
& +  \{-\alpha_1 \alpha_2   \alpha_3 \}     \nonumber 
\end{align}
for $  \alpha_1,  \alpha_2,  \alpha_3 \in  \mathbb C \mathbb P ^1$
. (The indices i are taken modulo 3 ).\\
Next he proves that we have the functional equation in 22 terms 
$$ (\ast \ast ) \qquad \qquad {\cal L}_3 ( R_3( a,b,c))=0  \quad
\quad  \qquad a,b,c \in \mathbb C $$  
Then he shows (part (a) of Theorem 1.10 in  \cite{gon} ) that 

 {\em ``the space of real continuous functions on
$\mathbb C \mathbb P ^1\setminus \{ 0,1, \infty \}$ that satisfy the
functional equation $(\ast \ast )$ is generated by the functions
${\cal L}_3(z) $ and $ {\cal L}_2(z).\log|z| $'' }.\\
(In fact, what he proves implies that this theorem is valid for measurable
 functions).\\

He had remarked before that, if we specialize this equation by setting
$ a=1, b=x, $  and  $ c=\frac{1-y}{1-x} $ , the equation $(\ast
\ast)$  simplifies and by using the inversion relation $ {\cal
  L}_3({x}^{-1})= {\cal L}_3(x) \: , \: x  \in  \mathbb C \mathbb P ^1
$, we obtain a homogeneous version
(i.e. without the right hand side $ {\sf R}_3(x,y)$) of equation~$ (SK)
$~.\\
This leads him to ask if this specialization characterizes 
the solutions of~$(\ast\ast)$.\\

The explicit determination of a basis of $ \underline{{\cal S}{ 
\stackrel{\cal O}{{}_\omega}}}   ({\cal SK} ) $    done in part 3.4. allows
us to give a positive answer to this question\,:
 we have this real semi-local characterization of  $  {\cal L}_3   $ :
\begin{theorem}
Let ${\cal G} :   ] -\infty,1 [ \:       \rightarrow \mathbb R$  be
a measurable function such that for $ 0<x<y<1$ we have 
\begin{align}
 2 \: {\cal G}\left(x\right)&+2 \: {\cal G}\left(y\right) - \: {\cal
  G}\left(\frac{x}{y}\right) +2 \: {\cal
  G}\left(\frac{1-y}{1-x}\right) +2 \: {\cal G}\left(\frac{x(1-y)}{y(1-x)}\right) 
- {\cal G}(xy)  \nonumber \\  
& +2 \: {\cal G}\left( \frac{x(1-y)}{x-1}      \right)
+2  \:  {\cal G}\left( \frac{y-1}{y(1-x)}\right) -{\cal
  G}\left(\frac{x(1-y)^2}{y(1-x)^2} \right)= 2 \l{3} (1)    \nonumber
\end{align}
 Then if we suppose ${\cal G}$ continuous at $0$, then there exists $\alpha\in \mathbb R$ such that $$ {\cal G}=\alpha
 \: {\cal L}_3+\frac{2}{9}(1-\alpha) \: \l{3}(1) $$
\end{theorem}

With our results of part 2 and 3.4, the proof 
is just a tedious exercise of linear algebra.\\
(It can be proved again by a suitable application of Abel's method).
It implies this result for $ \l{3} $ : 
\begin{coro}
Let $ F \! : \, ] -\infty,1 \: [ \: \rightarrow
\mathbb R $ be a measurable function such that for $ \, 
0<x<y<1 $, we have 
\begin{align}    
2  \, F(U_1(x,y)) & + \, 2 \, F(U_2(x,y))- \, F(U_3(x,y))      \nonumber \\
    & \qquad  +\,   2 \, F(U_4(x,y)) +\,   2 \, F(U_5(x,y)) -
  \, F(U_6(x,y)) \qquad    
\nonumber \\
 &   \qquad \qquad     + \, 2 \, F(U_7(x,y))+\,  2 \, F(U_8(x,y))-
   \,  F(U_9(x,y))= \, {\sf R}_3(x,y) \nonumber 
\end{align}
\begin{itemize}  \item If $F$ is continuous at 0 then there
  exists $ \alpha \in \mathbb R$ such that $$F=\l{3}+\alpha  \: ( {\cal
  L}_3-\frac{2}{9}\l{3}( 1) )$$
\item If $F$ is derivable at 0 then  $ F=\l{3} $ .
\end{itemize}
\end{coro}

\begin{flushleft} {\it Luc Pirio, \\
Equipe d'analyse complexe,\\
Institut de math{\'e}matiques de Jussieu, \\
175 rue du Chevaleret, 75013 Paris France\\
luclechat@hotmail.com\\
pirio@math.jussieu.fr}
\end{flushleft}
\end{document}